\documentclass[10pt]{article}
\usepackage{amsfonts,amsmath,mathrsfs,amssymb}
\usepackage{tikz}
\usepackage{verbatim}
\usepackage{array}
\usepackage[hang]{subfigure}
\usepackage{tabularx}

\newenvironment{myindentpar}[1]%
{\begin{list}{}%
         {\setlength{\leftmargin}{#1}}%
         \item[]%
}
{\end{list}}

\newtheorem{theorem}{Theorem}[section]
\newtheorem{lemma}[theorem]{Lemma}

\newtheorem{prop}[theorem]{Proposition}
\newtheorem{defin}{Definition}[section]

\newtheorem{cor}[theorem]{Corollary}

\def\R{ \mathbb{R}}

\def\proof{\noindent {\bf Proof.  }}

\def\eop{\hfill $\square$ \medskip}

\def\complaint#1{}
\def\withcomplaints{
\newcounter{mycomplaints}
\def\complaint##1{\refstepcounter{mycomplaints}%
\ifhmode%
\unskip%
{\dimen1=\baselineskip \divide\dimen1 by 2 %
\raise\dimen1\llap{\tiny -\themycomplaints-}}\fi%
\marginpar{\tiny [\themycomplaints]: \textcolor{revision} {##1}}}%
}

\definecolor{revision}{rgb}{1,0,0}
\definecolor{adnan}{rgb}{0,0,1}

\def\blocktriangular{block-triangular }

\def\dassur{$d$-Assur }
\def\tassur{$2$-Assur }
\def\thassur{$3$-Assur }
\def\sdassur{strongly $d$-Assur }
\def\sthassur{strongly $3$-Assur }
\def\dassurp{$d$-Assur}

\def\thassurp{$3$-Assur}
\def\sdassurp{strongly $d$-Assur}
\def\sthassurp{strongly $3$-Assur}

\makeatletter
\def\Ddots{\mathinner{\mkern1mu\raise\p@
\vbox{\kern7\p@\hbox{.}}\mkern2mu
\raise4\p@\hbox{.}\mkern2mu\raise7\p@\hbox{.}\mkern1mu}}
\makeatother

\begin{document}

\title {Directed Graphs, Decompositions, and Spatial Linkages}

\author{Offer Shai \thanks{Supported by a grant from ISF (Israel)} \\ Faculty of Engineering, Tel-Aviv University. shai@eng.tau.ac.Il
\\Adnan Sljoka\thanks{Supported in part under a grant from NSERC (Canada).}\\
 Department of Mathematics and Statistics\\ York University, 4700 Keele Street\\ Toronto, ON M3J1P3, Canada
and\\ Walter Whiteley\footnote{Supported by a grant from NSERC (Canada).}\\
 Department of Mathematics and Statistics\\ York University, 4700 Keele Street\\ Toronto, ON M3J1P3, Canada
}

\maketitle

\begin{abstract}
The decomposition of a linkage into minimal
components is a central tool of analysis and synthesis of
linkages.   In this paper we prove that every pinned $d$-isostatic
(minimally rigid) graph (grounded linkage) has a unique
decomposition into minimal strongly connected components (in the
sense of directed graphs), or equivalently into minimal pinned isostatic graphs,  which we call \dassur  graphs.
We also study key properties of motions induced by removing an
edge in a \dassur graph - defining a stronger sub-class of
\sdassur graphs by the property that all inner vertices go into
motion, for each removed edge. The \sthassur graphs are the
central building blocks for  kinematic linkages in $3$-space and
the 
 \thassur  graphs are components in the analysis of
built linkages. The \dassur graphs share a number of key
combinatorial and geometric properties with the \tassur graphs,
including an associated lower \blocktriangular decomposition of
the pinned rigidity matrix which provides modular information for
extending the motion induced by inserting one driver in a bottom
Assur linkage to the joints of the entire linkage.  We also
highlight some problems in combinatorial rigidity in higher
dimensions ($d$ $\geq$ 3) which cause the distinction between
\dassur and \sdassur which did not occur in the plane.
\end{abstract}


\section{Introduction}
\begin{figure}[h!]
\centering
   \subfigure[] { \includegraphics[width=.4\textwidth]{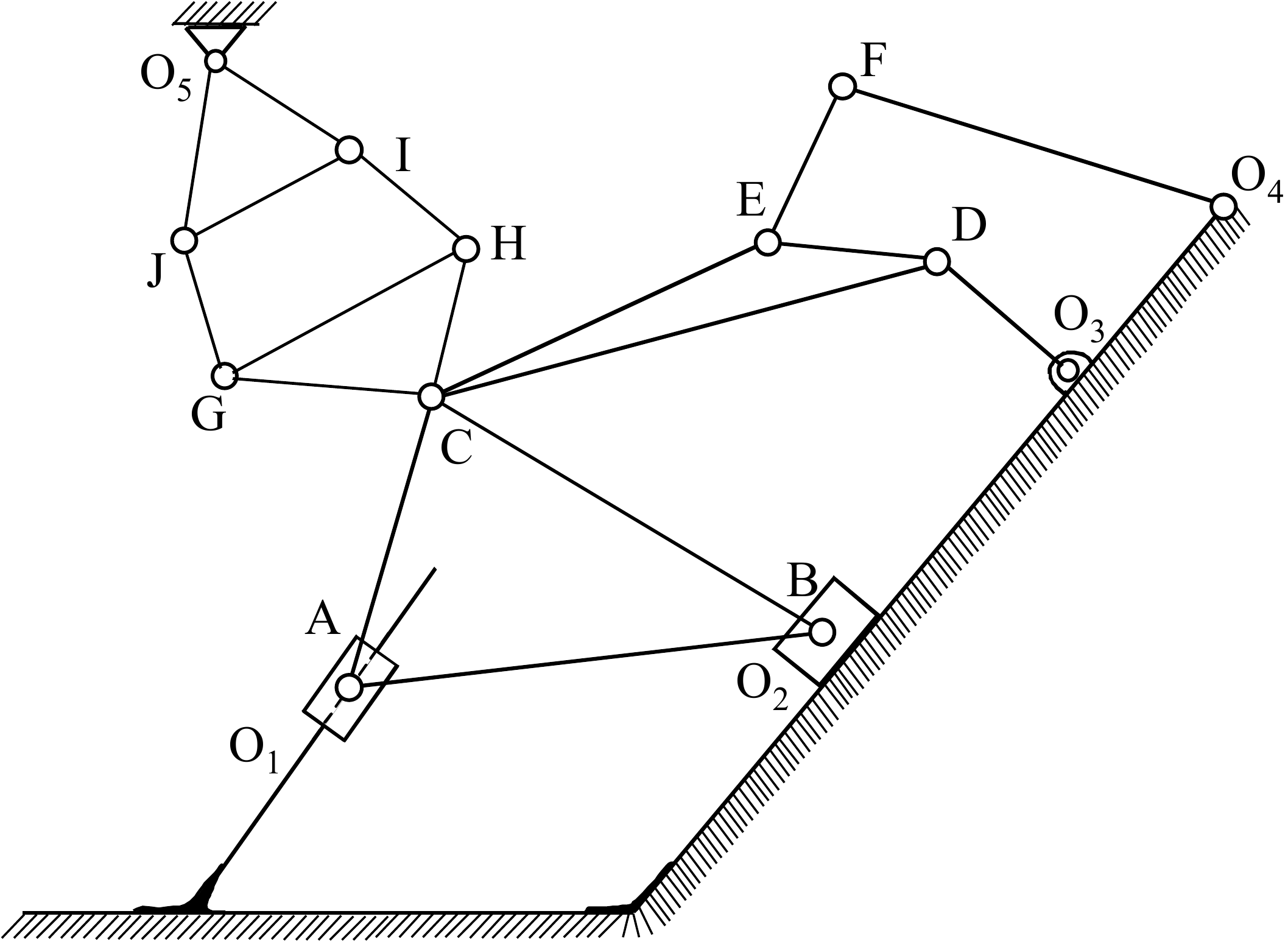}} \quad
    \subfigure[] { \includegraphics[width=.4\textwidth]{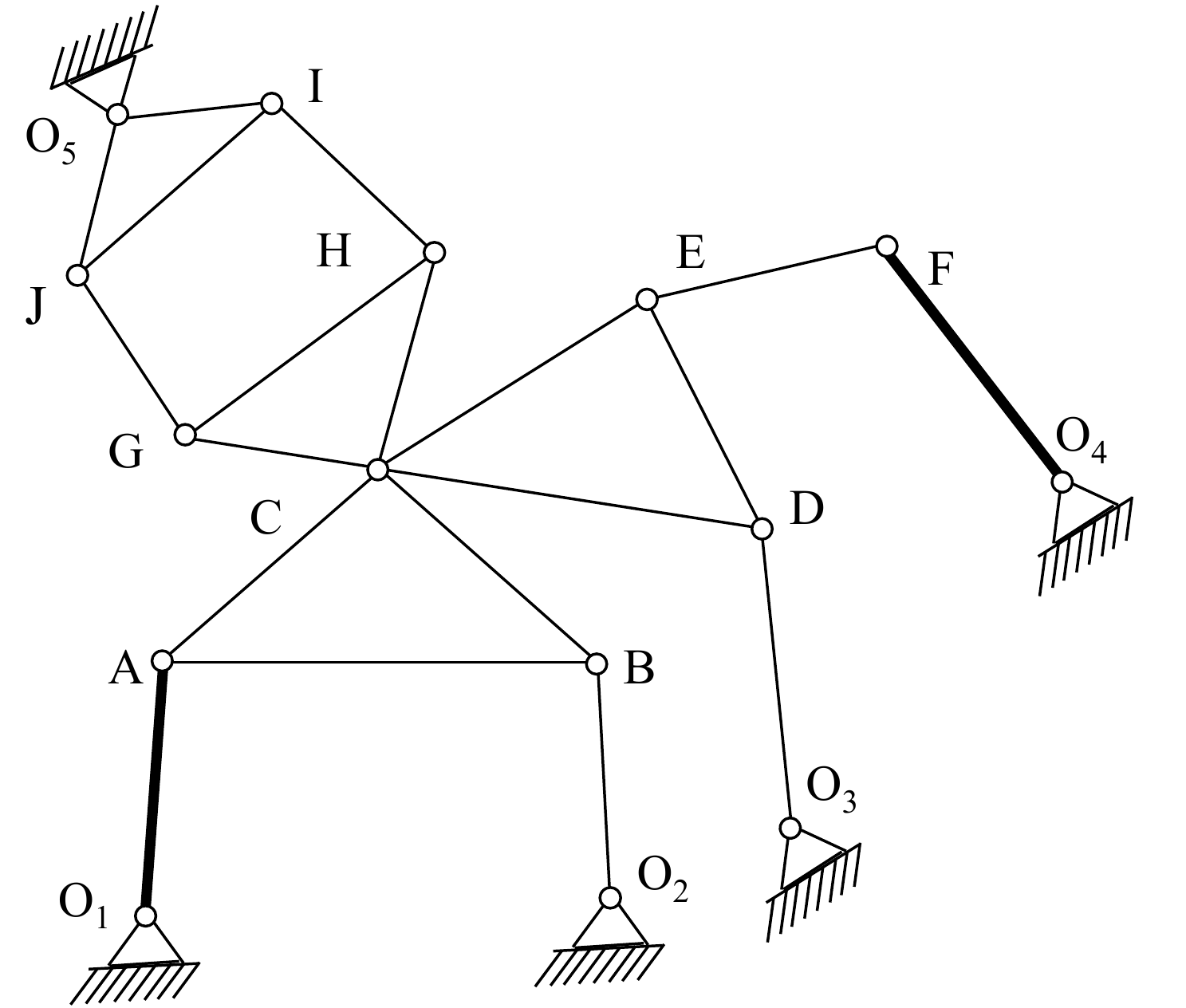}} \quad
   \subfigure[] { \includegraphics[width=.4\textwidth]{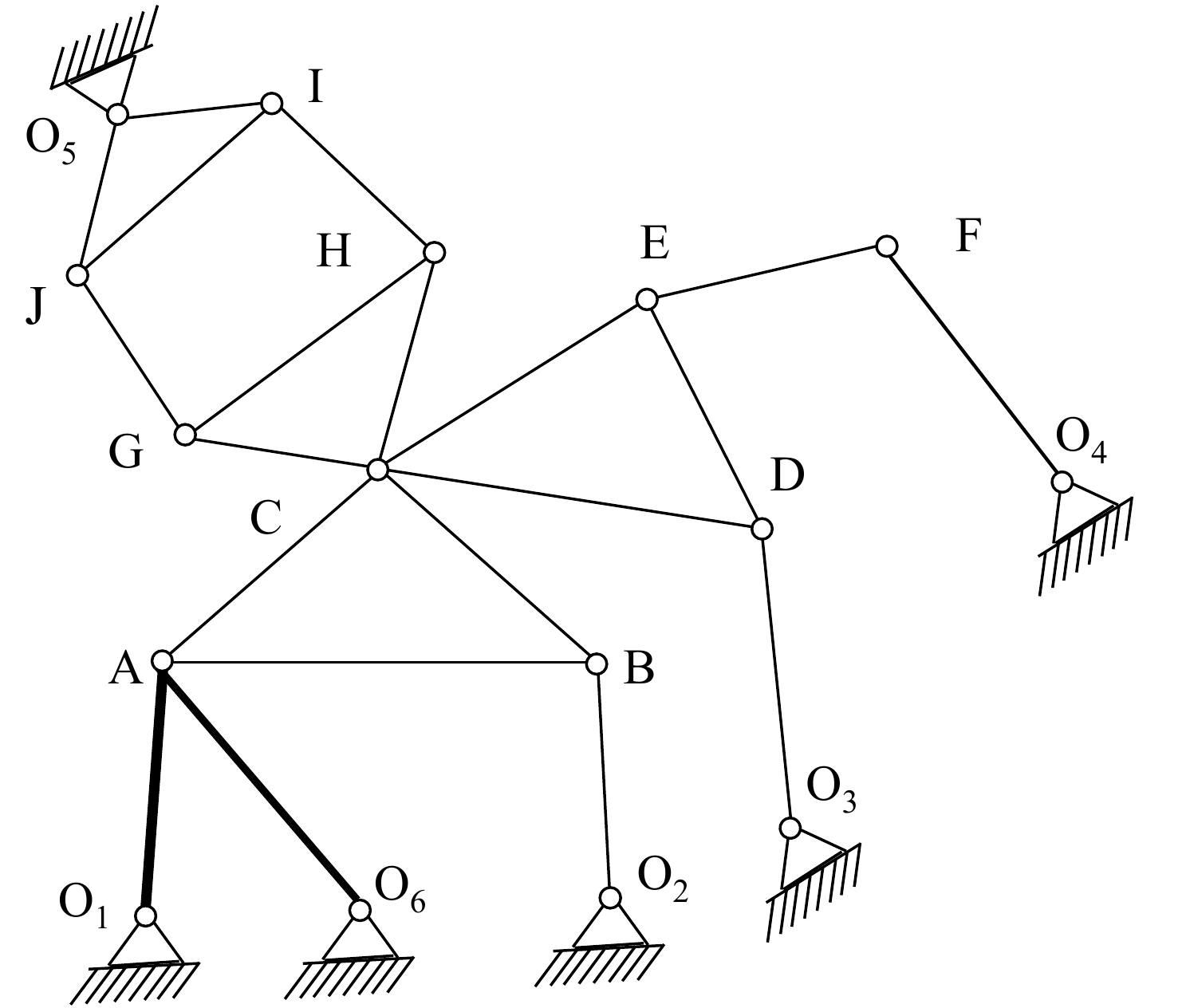}}\quad
    \subfigure[] { \includegraphics[width=.4\textwidth]{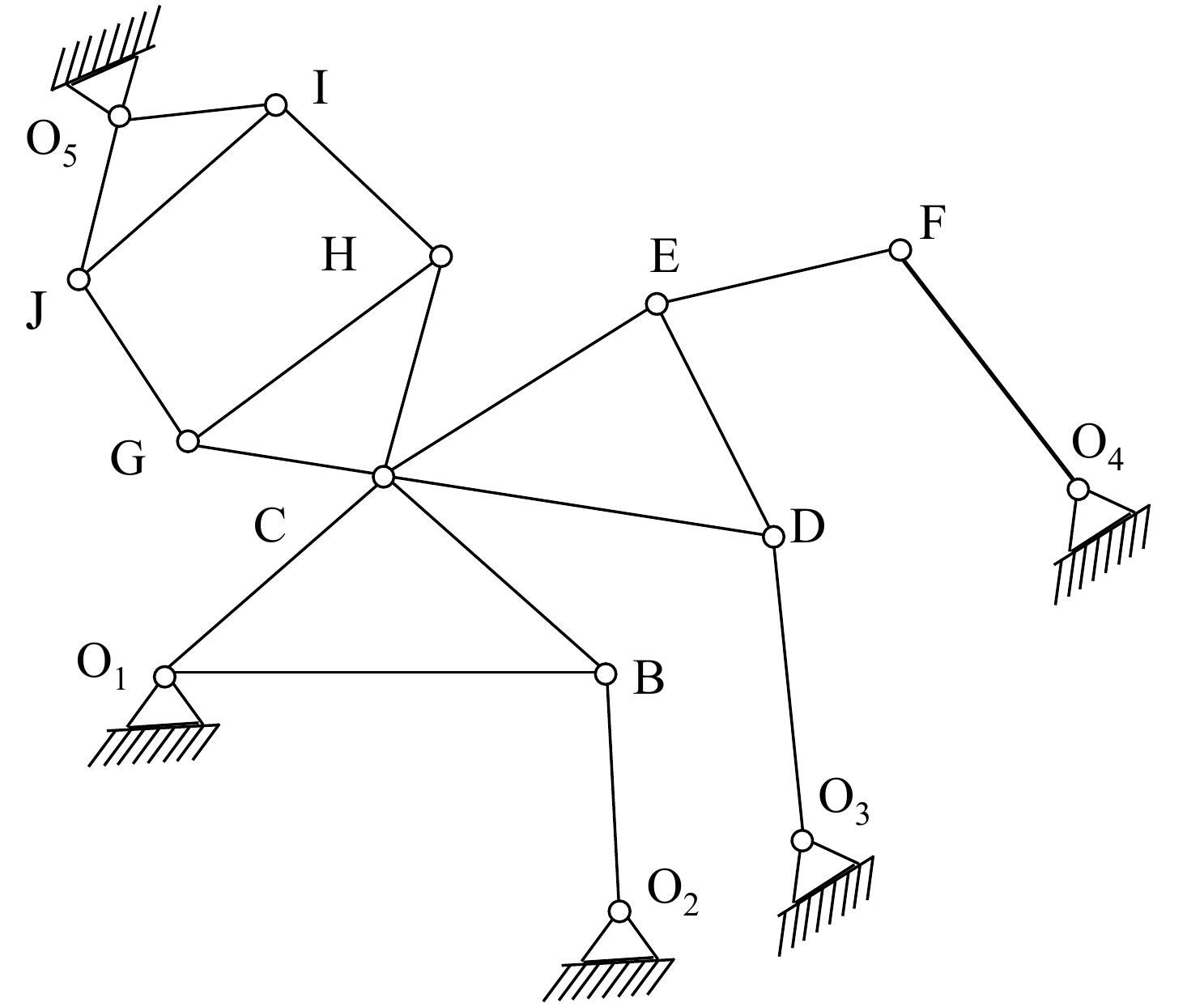}}\quad
    \subfigure[] { \includegraphics[width=.4\textwidth]{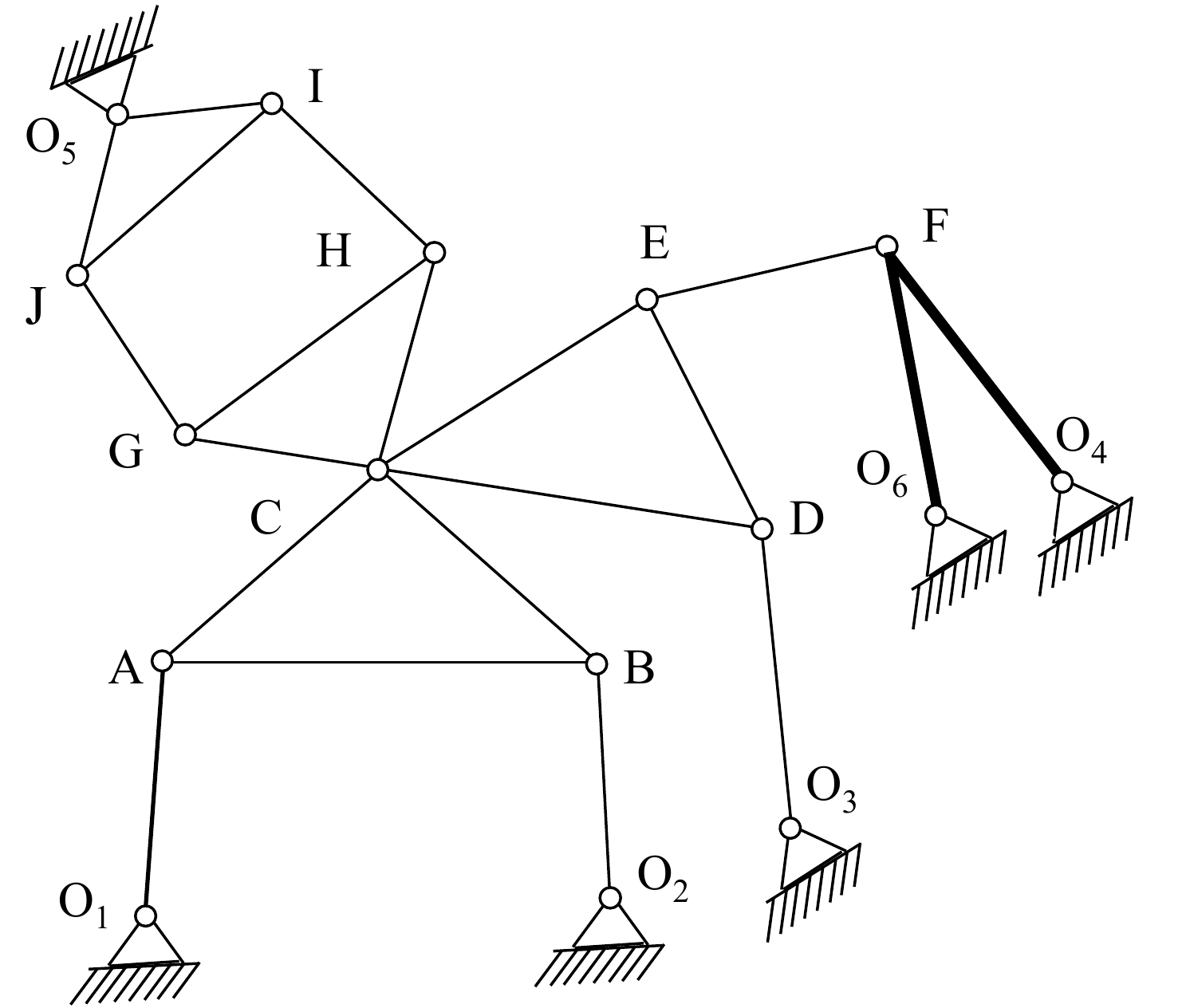}}\quad
    \caption{A plane linkage (a) is transformed into a flexible pinned framework (b). Assuming the chosen driver is the link from A to the ground, with one added bar (c) or by pinning the end of the driver to the ground (d) this becomes an isostatic pinned
    framework. If we specify the link from F to the ground as a designated driver, with one added bar, this results in an alternate pinned isostatic framework
    (e).}
    \label{fig:Linkage}
    \end{figure}

The decomposition of a system of constraints into small basic
components is an important tool of design and analysis.
In particular, the decomposition of a mechanical engineering  linkage into
minimal mechanical components is a central tool of analysis and
synthesis of mechanisms \cite{Assur,Mitsi,SSWI,SSWII}.  Figure~\ref{fig:Linkage} illustrates the analysis: the initial plane linkage (Figure~\ref{fig:Linkage}(a)),
is transformed into a flexible pinned framework
(Figure~\ref{fig:Linkage}(b));  designating one of the links
to be a driver, then adding an extra bar, or  pinning the end of
the driver to the ground, this becomes an isostatic pinned
framework (Figure~\ref{fig:Linkage}(c,d)) (see \S4).

The focus of
this paper is the study of specific decompositions of such
associated pinned isostatic frameworks  in dimensions 2 and 3.
However, the theorems and definitions will be generalized whenever
possible to all dimensions $d\geq 2$. We emphasize that this
decomposition goes further into the framework than the
decomposition of a flexible framework such as (b) into rigid
components.

For plane linkages, this decomposition has been an important tool
in the mechanical engineering literature over some decades \cite{Assur,Mitsi}. Our
central goal with this paper is to extend the previous
decomposition for plane pinned isostatic frameworks \cite{SSWI,SSWII} to pinned
isostatic frameworks in $d$-space
\cite{mechanicalpebblegame}. In developing this
extension, we present some key additional properties of the Assur
decomposition in the plane (and in $d$-space)  in terms of lower
\blocktriangular decompositions of the associated pinned rigidity
matrix (\S3). In the analysis, we also draw some new connections
to the theory of strongly connected decompositions of directed
graphs and their associated algorithms.

Section 2 highlights a unique decomposition of constraint graphs,
using only properties of the graphs of the constraints captured in
directions on the edges.  A key first step is generating a directed
graph for the constraints, which in \S3 will be directed towards
the `ground' of the linkage. This directed graph is then
decomposed into strongly connected components (components in
directed cycles in the directed graph), using a standard
combinatorial result which can be implemented using various
algorithms, such as Tarjan \cite{Tarjan}. Overall, the strongly
connected decomposition is presented as an acylic graph with
condensed nodes for the strongly connected components (Figure
\ref{fig:strongdecomp}).  In this decomposition, the strongly
connected components can be recognized visually as separated in the original
directed graph by directed cut-sets.

As a useful invariant property of this decomposition, we note that two equivalent orientations
of a given (multi)-graph (orientations with an assigned out-degree
for each vertex) will produce the same strongly connected
components, in the strongly connected decomposition.  This is
relevant to techniques in \S3 in this paper and to the
algorithms such as the pebble game, which include some choices of
orientations for edges -   choices that do not alter the
decomposition.

In Section 3, we show that for an appropriate orientation (a $d$-directed orientation ) this strongly connected graph
decomposition gives a decomposition of a pinned isostatic graph which coincides with several other key decompositions of a pinned isostatic graph: (i)
with a block triangular decomposition of the pinned rigidity matrix (\S3.2) and (ii) an associated decomposition into minimal pinned isostatic graphs (\S3.4).  For pinned linkages in the plane, this shared decomposition coincides with the $2$-Assur decomposition in
\cite{SSWI} .

Here we extend these decompositions to all higher dimensions $d\geq 2$.
The minimal pinned isostatic
 graphs in this decomposition are called \dassurp, and coincide with the plane $2$-Assur components of the previous work.
 This connection to the pinned rigidity matrix also provides a way to generate a
directed graph from the original undirected graph, so that we can
apply the strongly directed decomposition of \S2, if directions
were not already supplied by other analysis. This explicit
connection to lower triangular matrices is new in the plane, as
well.

In Section 4, we consider extensions to dimension $d$ of a further key property of
$2$-Assur graphs: removal of any single edge, at a generic
configuration, gives a non-trivial motion at all inner (unpinned)
vertices.   This property represents the desirable property that a
driver replacing this edge causes all parts of the mechanism to be
in motion. While this property follows for minimal pinned
isostatic graphs in the plane, it is a stronger property in higher
dimensions.  We use this added property to define the  \sdassur
graphs in all dimensions, as a restricted subclass of
\dassur graphs. This distinction between \dassur graphs
and \sdassur graphs adds another view on the complexity of
obtaining a combinatorial characterization of generic rigidity
3-space or higher. We will offer  examples which illustrate this
distinction.

This capacity to break the overall analysis of a pinned framework down into smaller
pieces, and recombine them efficiently, is the central contribution of these
decompositions to the analysis and synthesis of mechanical
linkages \cite{Mitsi,Norton,Shai}.  We conjecture that other
special geometric properties, such as those studied in
\cite{SSWII}, are also inherited by these \sdassur graphs in higher dimensions.

In \S5.1, we mention some further directions and extension of these
techniques and decompositions to the alternate body-bar frameworks
and related mechanisms.  In keeping with the general work on
rigidity in higher dimensions, these alternative structures have
good theories with \dassur again coinciding with
\sdassurp, and good fast algorithms for these structures.  Since
many built 3-D linkages have this modified structure, and these have
important applications to  fields such as mechanical
engineering and robotics, Assur body-bar frameworks and
decomposition algorithms are currently being actively explored.

We note that the strongly connected decomposition of a directed
graph is included in the code of current Computer Algebra Systems
such as Maple, Mathematica and SAGE. However, we alert the reader
that the lower \blocktriangular  decomposition for a matrix
implemented in these CAS systems is weaker than the decomposition
we present here for constraint graphs.  We comment further on this
in \S\ref{subsec:directiondmatrix}.

This general focus on decompositions of constraints into
irreducible components resembles other work done in electronics.
For example, Kron's diakoptics \cite{Kron}, describes a method
that decomposes any electrical network into sub-networks which can
be solved independently and then joined back obtaining the
solution of the whole network.  The overall goal of simplifying
analysis (and synthesis) through decomposition into pieces which
can be analyzed and then recombined through known ÔboundaryÕ
connections is part of a wide array of methods in systems
engineering \cite{Bowden}.

The directed graph decomposition presented in this
paper connects to general work on decomposition of CAD systems
\cite{OwenI,OwenII,Zhang}. In the latter systems the problem of
decomposing a cluster configuration into a sequence of clusters
\cite{HLSI,HLSII} corresponds to the decomposition of isostatic
graphs into \dassur graphs.  Some of that work also uses the
strongly connected decomposition highlighted
here 
as well as other related methods.

\noindent {\bf Aknowledgment} The authors thank the reviewer for a number of valuable
comments which were of great assistance to us.

\section {Decomposition of Pinned Directed  Graphs}\label{sec:directed}
In combinatorial rigidity theory, the rigidity of a given
framework is a property of an underlying undirected graph.
However, directed graphs come up in several important applications and algorithms
in rigidity theory.  For example, in the control theory of formations of autonomous agents,
the constraints on distances between agents are represented as a
directed graph \cite{Erin, Hendrickx}.
The output of the fast pebble game algorithm \cite{LeeStreinu,
mechanicalpebblegame, Sljoka} as part of the algorithmic
verification of critical counts and, implicitly, for decomposing
rigid and flexible regions in the framework graph is also represented
by a directed graph. Directed graphs also appear in t
mechanical engineering practice when synthesizing and analyzing
linkages.
We will use the directed graphs and their decompositions given in
this section  to obtain the \dassur decomposition of a
pinned isostatic framework in $d$-space in \S3.

In this section we state some basic definitions and background
from the theory of directed graphs
and present the decomposition of a directed pinned graph into
strongly connected components and develop a few simple extensions
to confirm the invariance of the decomposition under some natural
variations of the directed graph.

\subsection{Strongly connected component decomposition}

A \emph{graph} $G = (V, E)$ has a vertex set $V$ = $\{1,
2,...,n\}$ and an edge set $E$, where $E$ is a collection of
unordered pairs of \emph{vertices} called the \emph{edges} of the
graph.
We define a \emph{direction assignment} $\overrightarrow{G}$ to
graph $G = (V, E)$ as a pair of maps init:$E \rightarrow V$ and
ter:$E \rightarrow V$ assigning to every edge $e$ an \emph{initial
vertex} init($e$) and a \emph{terminal vertex} ter($e$). The edge
$e$ is said to be \emph{directed out} of init($e$) and {\it into}
ter($e$). 
We refer to $\overrightarrow{G}$ as a \emph{directed graph}
associated with $G$. A directed graph may have multiple directed
edges between the same two vertices, $v$ and $w$. We will call
such graphs directed \emph{multi-graphs}.

\medskip
\noindent {\bf Remark.}  For simplicity most definitions and
proofs are presented in terms of the vocabulary of simple graphs.
However, every proof in this paper will also apply to
multi-graphs. The functions init:$E \rightarrow V$ and
ter:$E\rightarrow V$ work well in this multi-graph setting, where
an edge is no longer uniquely represented by an ordered pair.
\eop

A \emph{cycle} of a graph $G$ is a subset of the edge set of G
which forms a path such that the start vertex and end vertex are
the same. A \emph{directed cycle} is an oriented cycle such that
all directed edges are oriented in the same direction along the cycle. A directed
graph $\overrightarrow{G}$ is \emph{acyclic} if it does not
contain any directed cycle. The \emph{degree} or \emph{valence} of
a vertex $v$ is the number of edges
that have $v$ as an endvertex. The \emph{out-degree} of a given
vertex in a directed graph is the number of edges directed out of
that vertex.
A vertex which has out-degree 0 is called a \emph{sink}. Sink
vertices will be \emph{pinned vertices} and directed graphs that have some
pinned vertices will be called \emph{pinned graphs}.

\begin{figure}[h!]
\centering
   \subfigure[] { \includegraphics[width=.4\textwidth]{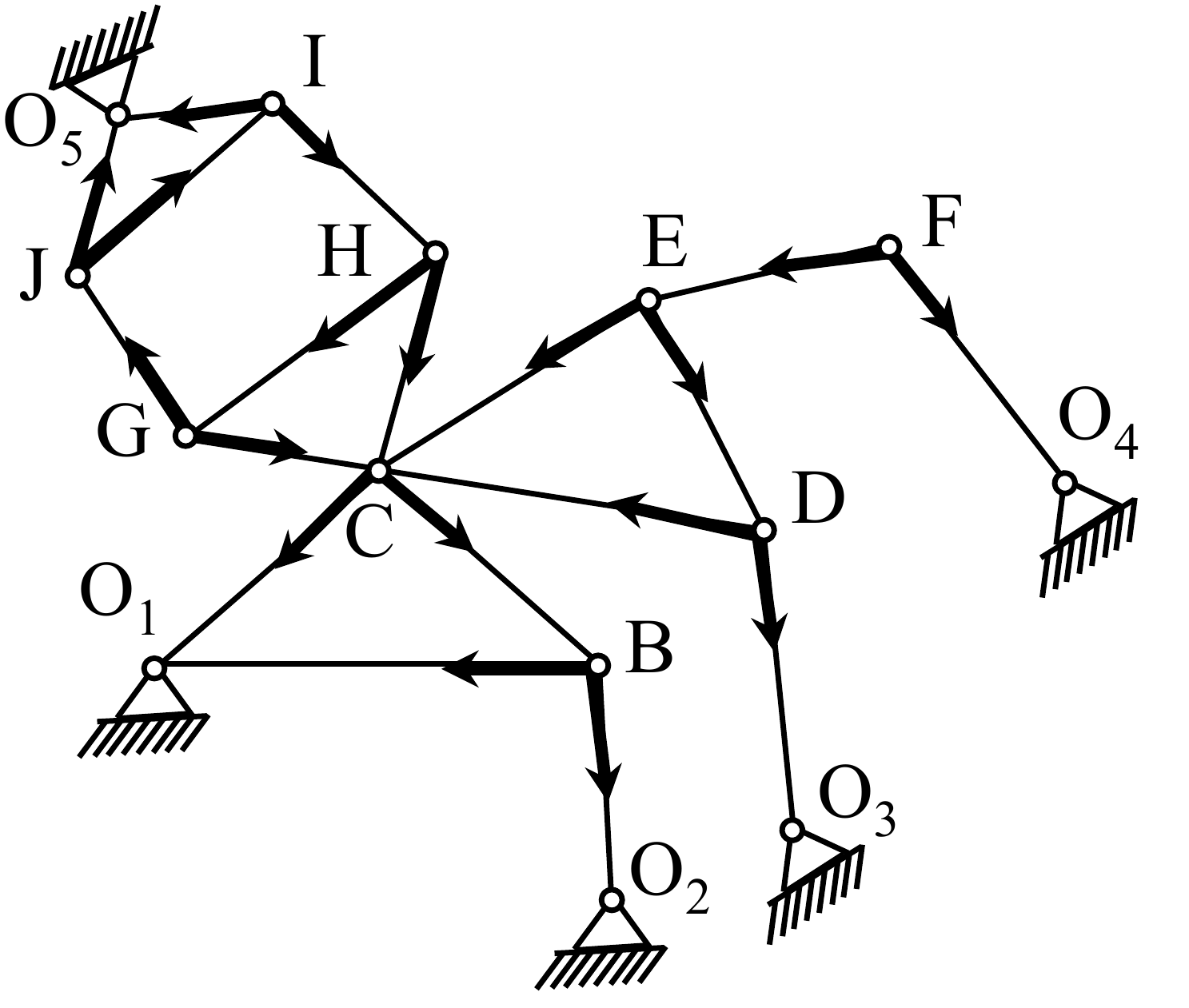}}
   \subfigure[] { \includegraphics[width=.34\textwidth]{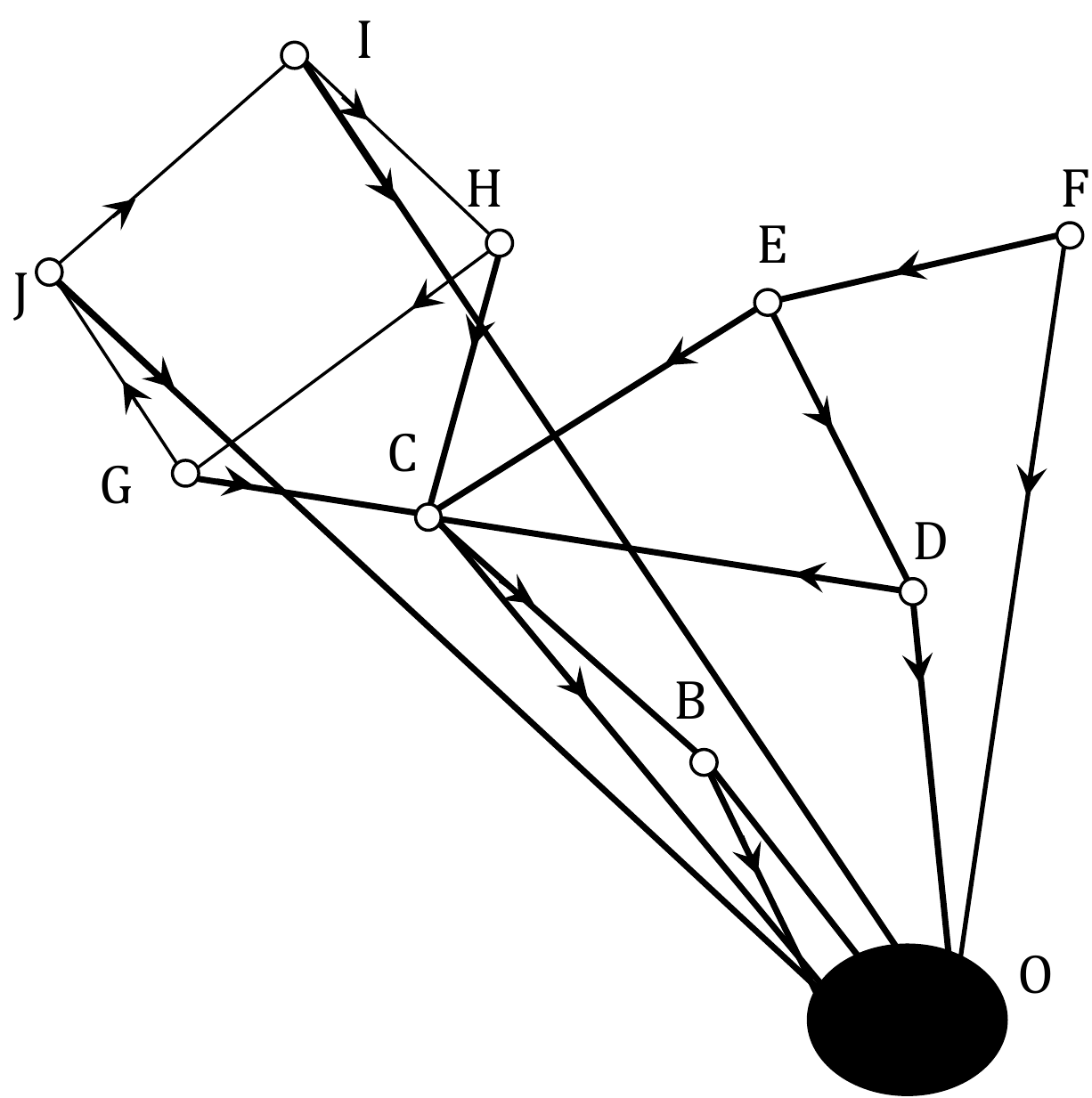}}
 \subfigure[] { \includegraphics[width=.36\textwidth]{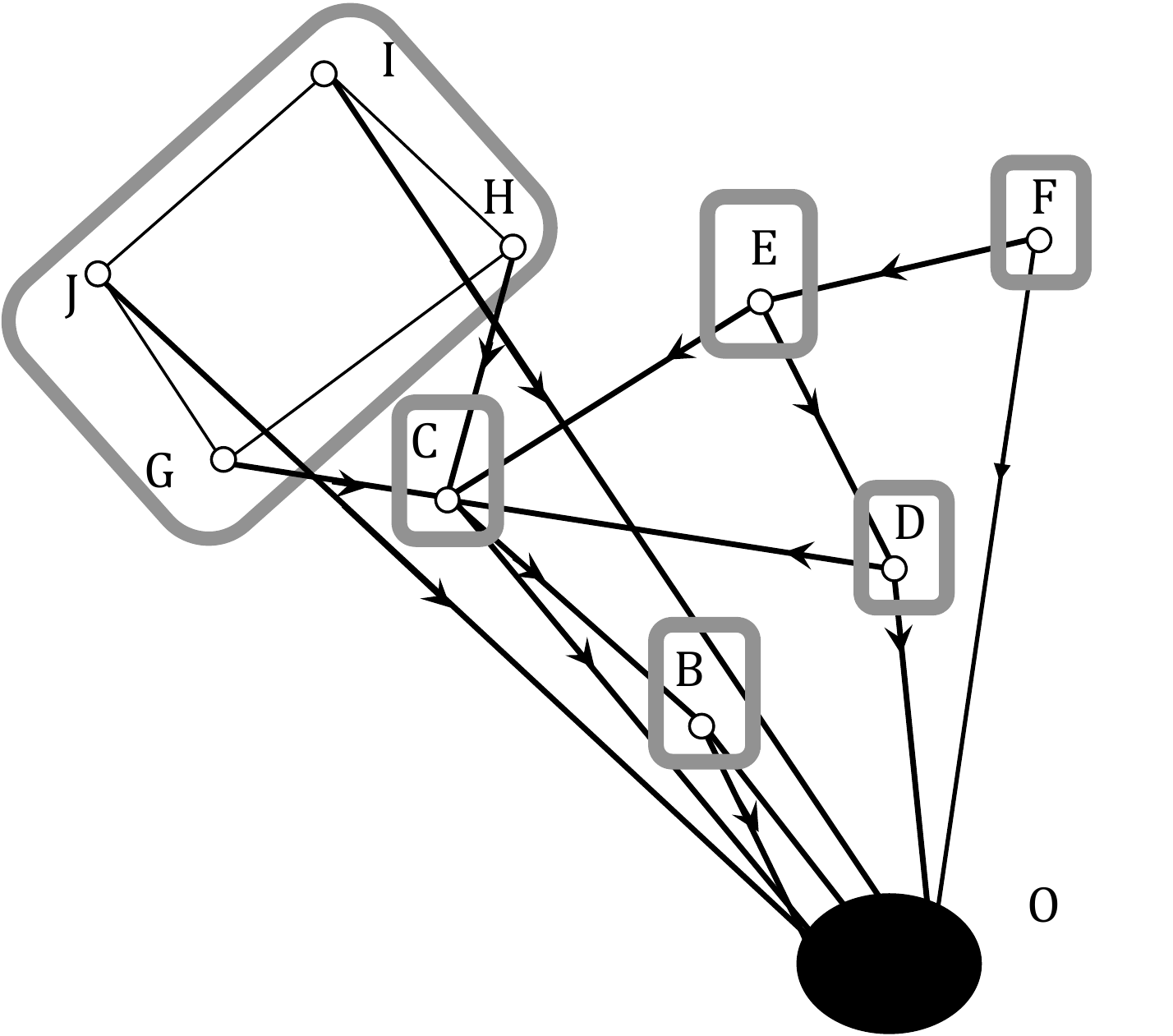}} \
   \subfigure[] { \includegraphics[width=.24\textwidth]{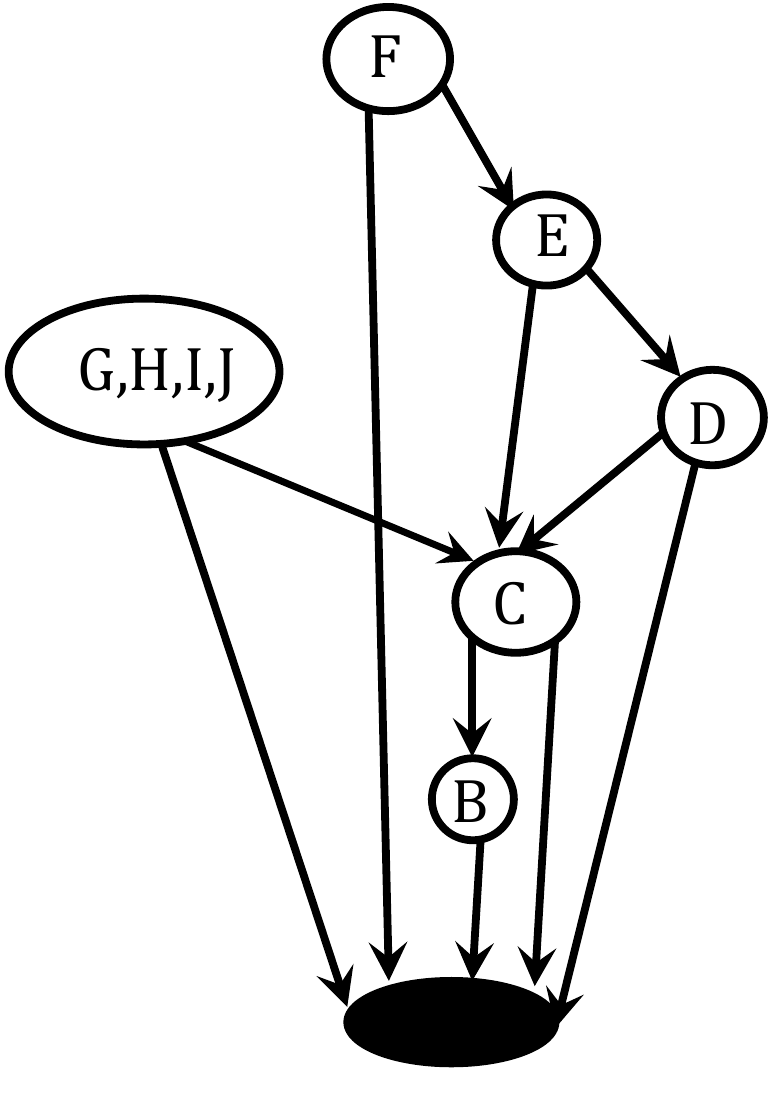}}
    \caption{Decomposition of a directed pinned graph: The directed pinned graph (a) has the pinned vertices condensed to the ground (sink) (b), with the corresponding strongly connected decomposition (c) and the partial order (d). }
    \label{fig:strongdecomp}
    \end{figure}

A directed graph is called \emph{strongly connected} if and only
if for any two vertices $i$ and $j$ in $\overrightarrow{G}$, there
is a directed path from $i$ to $j$ and from $j$ to $i$. The
\emph{strongly connected components} of a graph are its maximal
strongly connected subgraphs. It is maximal in the sense that a
subgraph cannot be enlarged to another strongly connected subgraph
by including additional vertices and its associated edges (see
Fig~\ref{fig:strongdecomp} (b), (c)).  One can determine the strongly connected components of a directed
graph using the $O$($|$E$|$) Tarjan's algorithm \cite{Tarjan}
which is implemented in several computer algebra packages, such as Maple,
Mathematica and SAGE. 

To illustrate the decomposition process of a pinned directed
graph, we first condense (ground) all pinned vertices into a
single ground vertex (sink)
(Figure~\ref{fig:strongdecomp}(a, b)). (More detailed definitions
of pinned graphs and grounding will be presented in the next
section.) After we have identified the strongly connected
components, as a simplification, we will ignore the orientation of
edges within the strongly connected components, and just keep the
orientation of edges between the components
(Figure~\ref{fig:strongdecomp} (c)). In the final step, we apply
condensation to each strongly connected component by contracting
it to a single vertex, obtaining an acyclic graph and a partial
order (Figure~\ref{fig:strongdecomp}(d)). We have schematized the
partial order, so that there are no multiple edges appearing
between any components.

So far, we were not concerned how we obtained the directed graphs.
In the next section, we will find directed graphs for our
linkages, and use those for our decomposition.  There are also
inductive constructions \cite{TayWhiteley, WChapter}, which are
currently being further investigated, and there are additional fast algorithmic techniques
for generating the desired directed graphs for linkages (e.g.
the pebble game) 
\cite{mechanicalpebblegame}.

\subsection{Equivalent orientations of a graph}
We now show that choices in orientations of edges which conserve
a fixed out-degree of each of the vertices do not alter the
decompositions. This is a useful observation for the proofs in the
next section and will assist the verification of algorithms for
generating the directed graphs and related decompositions
\cite{mechanicalpebblegame}. For example, all plays of the pebble
game algorithm with the same number of final pebbles at the
corresponding vertices will give the same strongly connected
decompositions \cite{LeeStreinu, Sljoka}. We thank Jack Snoeyink
for conversations which clarified these arguments.

\begin{defin}
Given a multi-graph $G$ and two direction assignments
${\overrightarrow{G^{1}}}$ and ${\overrightarrow{G^{2}}}$. We say
that ${\overrightarrow{G^{1}}}$ and ${\overrightarrow{G^{2}}}$ are
{\rm equivalent orientations} on G if the corresponding vertices
have the same out-degree.
\end{defin}

Such alternate orientations appear when comparing assigned
direction assignments by hand (common engineering practice) vs
those attained by algorithms. We confirm we obtain same decomposition  regardless how we
obtained the equivalent direction assignment.

\begin{lemma}\label{lemma:eqorientations} Given two equivalent orientations  ${\overrightarrow{G^{1}}}$ and ${\overrightarrow{G^{2}}}$ on $G$, then the two orientations differ by reversals on a set of directed cycles.
\end{lemma}

\begin{figure}[ht]
\centering
    \subfigure[] {\includegraphics [width=.25\textwidth]{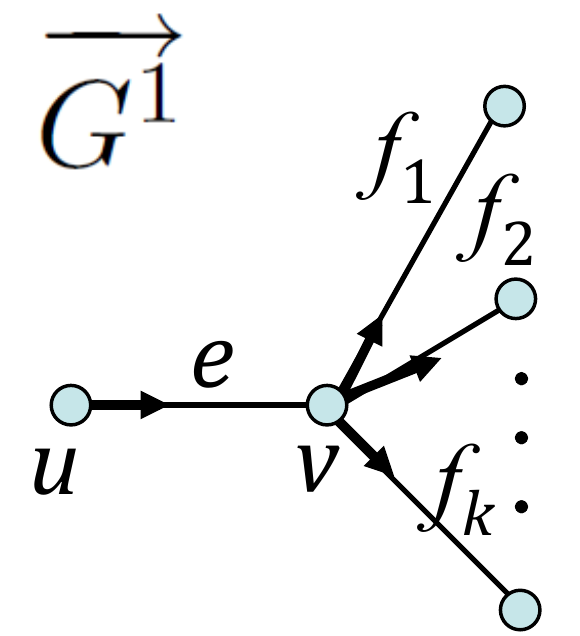}} \quad
    \subfigure[] { \includegraphics[width=.25\textwidth]{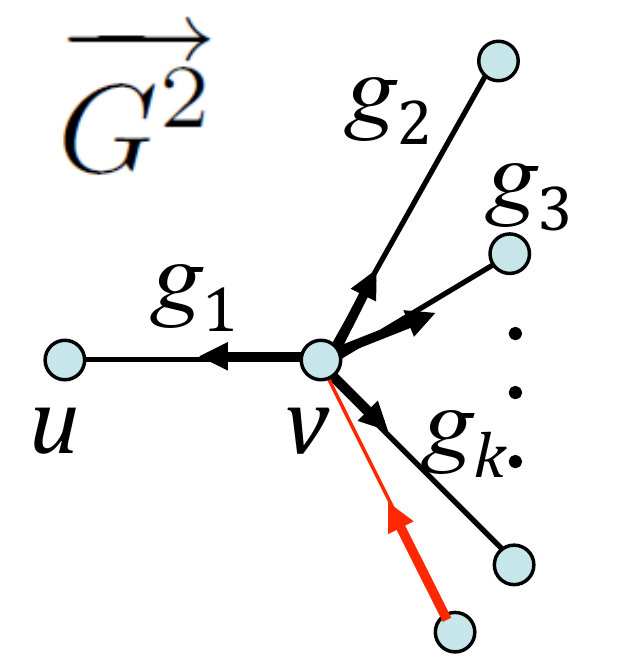}}\quad\quad
   \subfigure[] { \includegraphics[width=.24\textwidth]{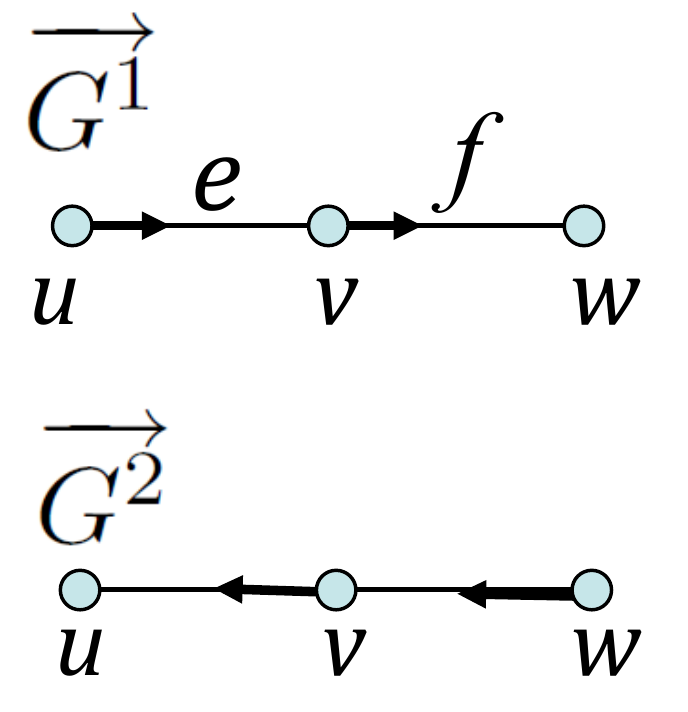}}
    \caption{Since vertex $v$ has same out-degree in ${{G^{1}}}$ and
in ${{G^{2}}}$,  there exists an edge in ${{G^{1}}}$ (a) that is
oppositely oriented in ${{G^{2}}}$ (b). We show such a selection
$f$ in (c,d).}
    \label{fig:pathreversalproof}
    \end{figure}

\proof Pick an edge $e$ = $(u, v)$ in ${\overrightarrow{G^{1}}}$
that is oppositely directed in ${\overrightarrow{G^{2}}}$. So, in
${\overrightarrow{G^{1}}}$ edge $e$ is incoming at vertex $v$.
Assume there are $k$ outgoing edges at $v$ in
${\overrightarrow{G^{1}}}$ Because $v$ has same out-degree $k$ in
${\overrightarrow{G^{2}}}$, and this edge is reversed to be
outgoing in ${\overrightarrow{G^{2}}}$ there must exist an
out-going edge from $v$ in ${\overrightarrow{G^{1}}}$, say $f$ =
$(v, w)$,  that is oppositely directed in
${\overrightarrow{G^{2}}}$
(Figure~\ref{fig:pathreversalproof}(a,b)).  We walk out of $v$
along $f$ in  ${\overrightarrow{G^{1}}}$
(Figure~\ref{fig:pathreversalproof}(c)).

As we enter a new vertex, we again apply this same argument.  This
identifies a directed path in ${\overrightarrow{G^{1}}}$ that is
oppositely directed in ${\overrightarrow{G^{2}}}$. As there is
only a finite collection of edges that have opposite direction, we
walk along this directed path until we come back to some vertex on
this path, identifying a directed cycle in
${\overrightarrow{G^{1}}}$ that has an opposite orientation in
${\overrightarrow{G^{2}}}$.

\begin{figure}[h!]

\centering
   {\includegraphics [width=.90\textwidth]{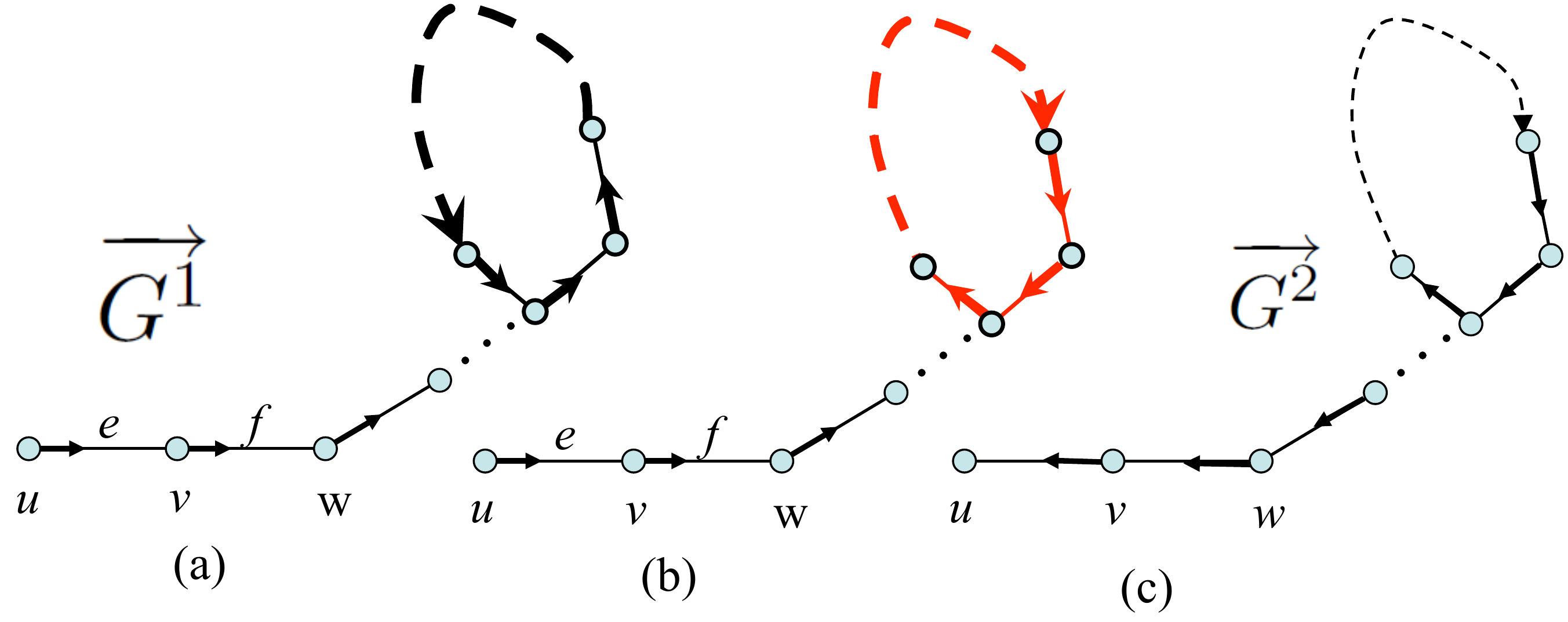}}
    \caption{Locating a directed path in $G^{1}$ that is
    oppositely oriented in $G^{2}$ (a), reversing a first cycle in
    $G^{1}$ (b), which now has the same orientation as in $G^{2}$
    (c).
}
    \label{fig:cyclereversalproof}
    \end{figure}

We reverse the orientation of such a cycle from the orientation in ${\overrightarrow{G^{1}}}$ towards the orientation in ${\overrightarrow{G^{2}}}$ This decreases the number
of edges in ${\overrightarrow{G^{1}}}$ that are oppositely
directed from ${\overrightarrow{G^{2}}}$ (see
Figure~\ref{fig:cyclereversalproof}).  We continue reversing
identified cycles, until all edges are directed following
${\overrightarrow{G^{2}}}$.      \eop

We end with a valuable corollary to this lemma:

\begin{cor}\label{lemma:samestronglycc}  Given two equivalent orientations ${\overrightarrow{G^{1}}}$ and ${\overrightarrow{G^{2}}}$,
the strongly connected components of the decompositions are the
same.
\end{cor}

\proof From Lemma~\ref{lemma:eqorientations},
${\overrightarrow{G^{1}}}$ and ${\overrightarrow{G^{2}}}$ differ
by reversal of set of directed cycles. As cycle reversals do not
change the strongly connected components the decompositions are
the same.\eop

\section{Decomposition of Pinned $d$-isostatic Graphs }  

Paper \cite{SSWI} described the 2-Assur decomposition of a pinned
generically isostatic graph in the plane.  Here we will nuance this decomposition with the connections to the directed graph decomposition of \S2,  and directly extend the central decomposition  to all dimensions.
This extended decomposition carries significant properties of the
$2$-Assur decomposition, with one key exception which we examine in \S4.
  Because
of applications in mechanical engineering, we are primarily
interested in the Assur decompositions in 2-space and 3-space. For
this reason all the figures will be either 2 or
3-dimensional.

Recall that the $2$-Assur graphs shared a suite of equivalent properties \cite{SSWI, SSWII}. We will begin with the extended definitions  in $d$-space, collectively defining the \dassur graphs (which coincide with Assur graphs for dimension $2$).  All the results in this section apply
to \dassur graphs. Accordingly, when we just speak of Assur graphs in
this section, it should be understood that we are referring to
\dassur graphs. For mathematical completeness, whenever
possible we state all the
theorems and definitions in $d$-dimensional space.

In \S4 we will explore one further property of \tassur graphs that does not always extend to higher dimensions - the response of inner vertices to removing one edge.

\subsection{Pinned $d$-isostatic graphs and pinned rigidity matrix}

\begin{figure}[htb]
\centering
 \subfigure[] {\includegraphics [width=.4\textwidth]{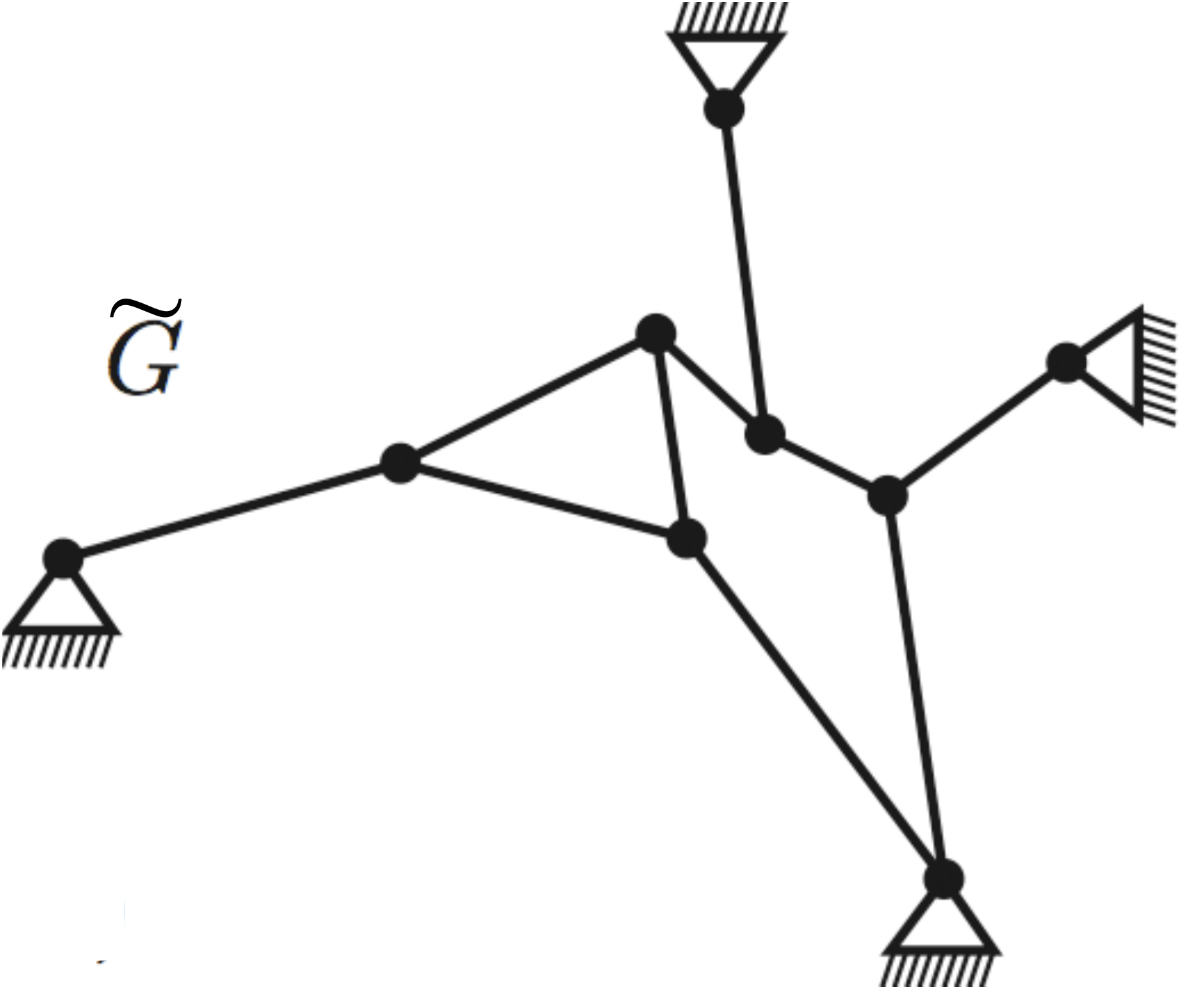}} \quad
   \subfigure[] { \includegraphics[width=.42\textwidth]{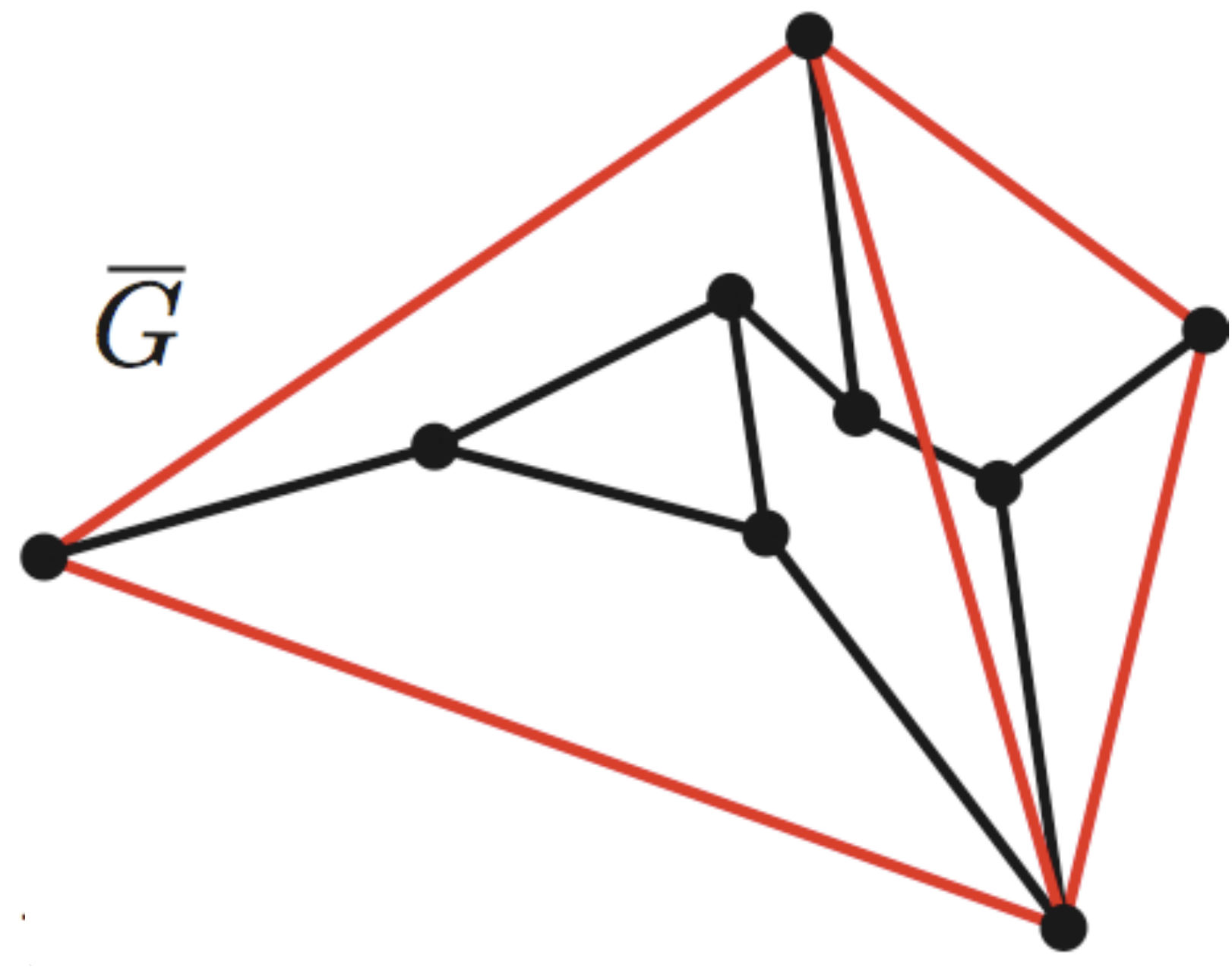}}
\caption{The framework (a) is pinned $2$-isostatic because
framework (b) is $2$-isostatic.\label{pinsub01}}
\end{figure}

In the remainder of the paper, we will focus on pinned frameworks
and the analysis of their associated graphs and directed graphs.  The rational for
this focus on pinned frameworks is as follows:

\vspace{-.3cm}
\begin{myindentpar}{.5cm}
 \noindent Given a framework associated with a linkage,  we are interested in its internal motions,
not the trivial ones.  Following the mechanical engineers, we pin
the framework
 by prescribing, for example in $2$-space, the coordinates of the endpoints of some
edges, or equivalently, by fixing the position of the vertices of
some rigid subgraph  (see Figure~\ref{pinsub01}).
\end{myindentpar}

We call these vertices with fixed positions {\em pinned}, the
other vertices are {\em inner}. (Inner vertices are sometimes called {\em
free} or {\em unpinned} in the literature on linkages.) Edges
between pinned vertices are irrelevant to the analysis of a pinned
framework. More formally, we denote a {\it pinned framework} as
$(\widetilde{G},p)= ((I,P;E),p)$, where $\widetilde{G}=(I,P;E)$ is
a pinned graph, $I$ is the set of inner vertices, $P$ is the set
of pinned vertices, $E$ is the set of edges, where each edge has
at least one endpoint in $I$, together with an assignment  $p$ of
points in $d$-space to the vertices of $\widetilde{G}$ (i.e.
$p$ is a fixed configuration (embedding) of $I\cup P$ into
$\mathbb{R}^{d}$).

In the previous papers \cite{SSWI, SSWII}, some of the rigidity analysis was routed through
grounding the pinned graphs through an
isostatic framework for the ground (Figure~\ref{pinsub01}).  Here we will give an
alternative version that directly analyzes the pinned rigidity
matrix, which only has columns for the inner vertices,
corresponding to the possible variations in their positions, if
the linkage is flexible with the pin fixed.

For a pinned framework $(\widetilde{G},p) = ((I,P;E),p)$ we define the
$|E|\times d|I|$ $d$-space {\it pinned rigidity matrix}, which
unlike the regular rigidity matrix, only has columns for the inner
vertices:
\begin{displaymath}
\mathbf{R}(\widetilde{G},p)\ =\ \bordermatrix{& & & & i & & & & j
& & & \cr
 & & \ddots& &  & & \vdots & &  &&\Ddots  &
\cr \{i,j\} & 0 & \ldots &  0 & (p_{i}-p_{j}) & 0 & \ldots & 0 &
(p_{j}-p_{i}) &  0 &  \ldots&  0 \cr & &  \vdots & &  & & \vdots &
&  &&  \vdots  & \cr \{i,k\} & 0 & \ldots &  0 & (p_{i}-p_{k}) & 0
& \ldots & 0 & 0 &  0 &  \ldots&  0 \cr & & \Ddots & &  & & \vdots
& &  & &\ddots  &
 }
\textrm{,}\end{displaymath} where $i,j \in I$ and $k\in P$. Note that this matrix has $d$ columns for each inner vertex.

The solutions $U$ to the equation: $\mathbf{R}(\widetilde{G},p)\times U^{tr}=0$ are called {\it infinitesimal motions} of the pinned framework.
A framework $(\widetilde{G},p)$ is {\em pinned $d$-rigid} if the only infinitesimal motion is the zero motion or equivalently, if
pinned rigidity matrix $\mathbf{R}$$(\widetilde{G},p)$ has full
rank $d|I|$. A framework is {\em pinned $d$-independent} if the
rows of $\mathbf{R}$$(\widetilde{G},p)$ are independent. A
framework $(\widetilde{G},p)$ is {\em pinned $d$-isostatic} if it
is both {\em pinned $d$-rigid} and {\em pinned $d$-independent}.
In particular, if the framework is pinned $d$-isostatic then
$|E|=d|I|$ and the pinned rigidity matrix is a square matrix.

If we vary the configuration $p$ over all of $\R^{d|I|+d|P|}$, then the pinned
rigidity matrix achieves some maximal rank, and this maximal rank
occurs for an open dense subset of $\R^{d|I|+d|P|}$ - the {\em
generic rank} of the {\em $d$-space pinned rigidity matrix} for the
graph.  In particular, for a pinned isostatic graph, the configurations that drop the  rank are captured by a non-zero
polynomial in variables for the vertices.  We conclude that if one
configuration $p$ achieves the full rank $d|I|$, then almost all
configurations (all points in this open dense subset) achieve this
rank, and we call all configurations in the open dense subset {\it
generic} or {\it regular}.

We say a pinned framework $(\widetilde{G},p)$ is {\it generic} if
the configuration  $p$ of the joints $p$ is generic. (A special subset of the generic configurations is those whose
coordinates satisfy no algebraic equations - also a dense set, but not an open set.)

\begin{theorem}\label{3PinEquivalenceThm}
Given a pinned graph $\widetilde{G} = (I,P;E)$, the following are
equivalent:

\begin{enumerate}
\item There exists a pinned $d$-isostatic realization $p$ of
$\widetilde{G}$ in $d$-space;

\item For all placements ${p}|_{P}$ of the pins $P$ in generic position in
$d$-space, and all generic positions of vertices in $I$ the
resulting pinned framework is pinned $d$-isostatic;

\item $\widetilde{G} = (I,P;E)$ is generically pinned $d$-rigid
and $|E|=d|V|$;

\item $\widetilde{G} = (I,P;E)$ is generically pinned
$d$-independent and $|E|=d|V|$.

\end{enumerate}
\end{theorem}

\proof  These statements are a simple translation of standard
results for isostatic frameworks to the pinned $d$-isostatic
setting \cite{WChapter}. \eop

Results in the plane, and other recent work suggests the {\it conjecture}  that (2) can be refined to require only that the pins are in general position within a $d$-hyperplane (a line in $\R^{2}$, a plane in $\R^{3}$).

We call any graph $\widetilde{G} = (I,P;E)$ satisfying the equivalent conditions of
Theorem~\ref{3PinEquivalenceThm} {\it pinned $d$-isostatic}.
A  {\em \dassur graph} is a \emph{minimal}
pinned $d$-isostatic graph.  By {minimal} we mean there is no
proper subgraph
which is also a pinned $d$-isostatic graph. In \S4, we will define a \sdassur
graph as a \dassur graph with the added property that
removal of any edge puts all inner vertices in motion. We will
look at the difference between the two types of graphs in more
detail in Section 4, but note no distinction between
$2$-Assur graphs and strongly $2$-Assur graphs (see also Section 4), so they will
always be called $2$-Assur.

\subsection{Directed graph decomposition of  the pinned rigidity matrix}\label{subsec:directiondmatrix}

We now use the $d$-space pinned rigidity matrix to
generate a special type of directed graph for any pinned graph in $d$-dimension
with $|E|=d|V|$.  We  then connect the corresponding directed graph decomposition from \S2 and  to a block decomposition of the pinned rigidity matrix, which will be central to the rest of the paper.

Matching the shape of the pinned rigidity matrix,
we develop {\em $d$-directed orientations of the pinned
graph}: one in which each inner vertex has out-degree $d$ and each
pinned vertex has out-degree $0$ - a sink in the directed graph.

\begin{prop} Every pinned  $d$-isostatic graph 
has a $d$-directed orientation, with all inner vertices of out-degree $d$ and all pinned vertices of out-degree $0$.  \label{prop:2orientation}
\end{prop}

\begin{figure}[htb]
\centering
  \subfigure[] {\includegraphics [width=.41\textwidth]{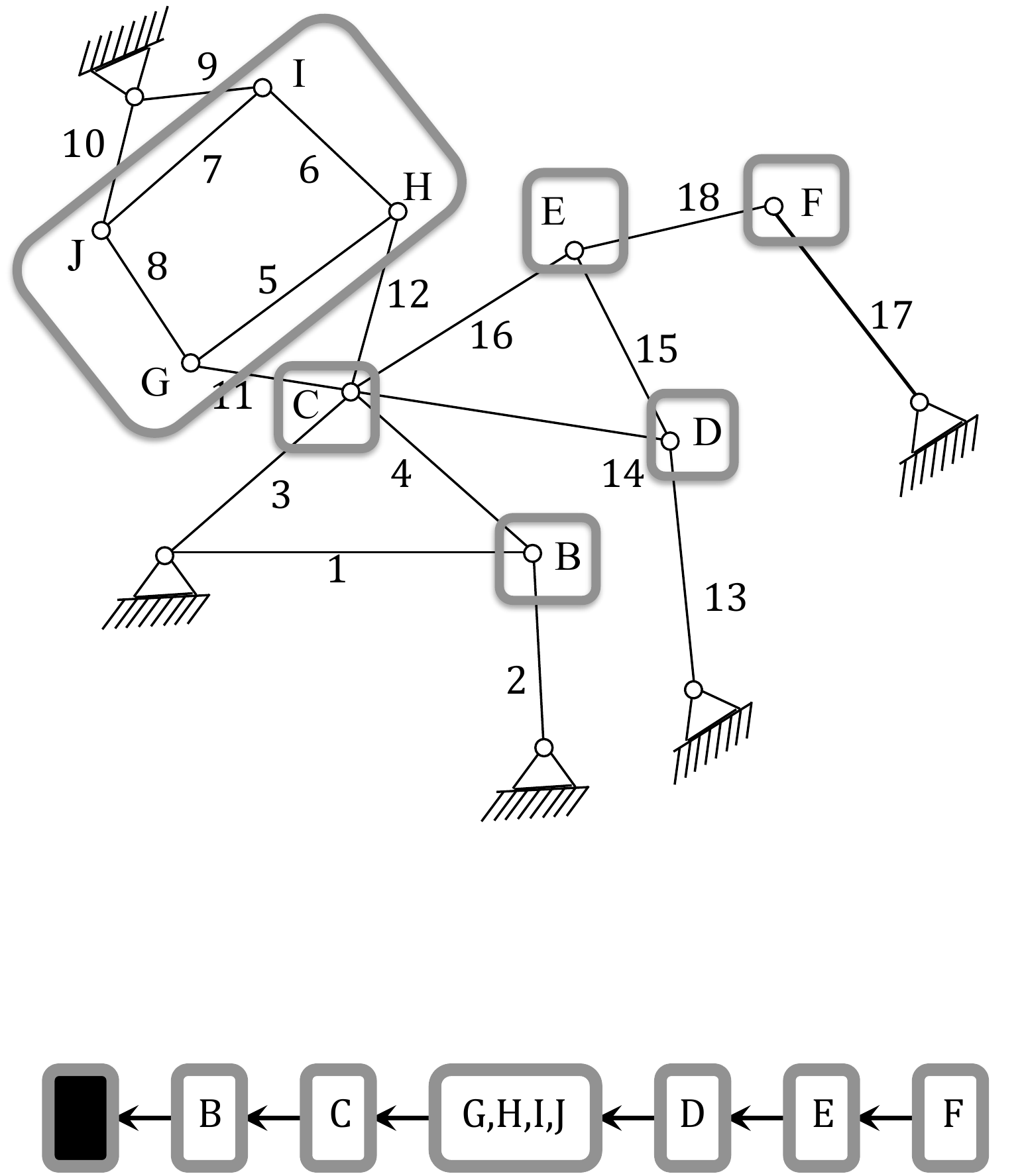} }\!\!\!\!\!\!
   \subfigure[] { \includegraphics[width=.58\textwidth]{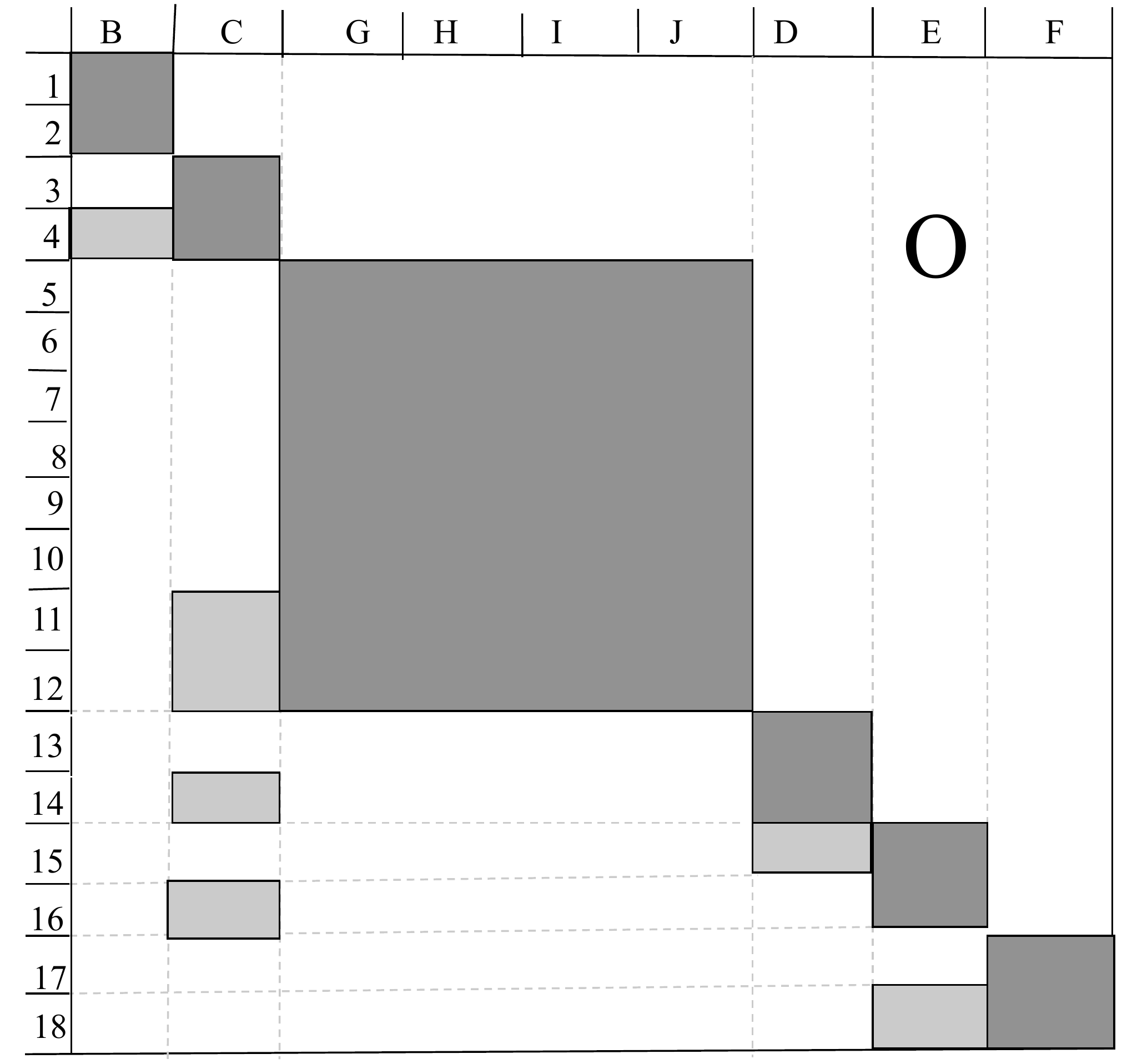} }
 \caption{ The pinned $2$-isostatic framework in (a), with the directed graph decomposition  of Figure~\ref{fig:strongdecomp} as a partial order,  extended to a linear order  generates a block triangular matrix in (b). \label{Block}}
\end{figure}

\proof Take the determinant of the $d|V|\times d|V|$ pinned
rigidity matrix.   We know the determinant of the square matrix is
non-zero if and only if the framework is infinitesimally rigid.
For this determinant to be non-zero, there must be a non-zero term
in the Laplace expansion of the determinant of this matrix, in
$d\times d$ blocks following the $d$ columns for each inner vertex.
Take any such non-zero term. This will associate $d$ rows (edges)
with each inner vertex $i$, and we direct these d edges out from
vertex $i$ (Figure~\ref{Block}).  This gives the desired
$d$-directed orientation. \eop

We can now apply the strongly connected graph decomposition from \S2 to
decompose the pinned $d$-isostatic graph, with the pinned vertices
identified in the directed graph decomposition.  Each strongly connected component is extended to include the outgoing edges from the component.   (In \S3.3, we show that these extended components are minimal
pinned $d$-isostatic graphs.) We make the connections through a
\blocktriangular  decomposition of the pinned $d$-rigidity matrix (Figure~\ref{Block}). In this decomposition, a permutation of the vertices and the edges of the pinned graph generates a lower block triangular matrix.

\begin{theorem} For a pinned $d$-isostatic graph with a pinned
$d$-directed orientation (all inner vertices of out-degree $d$ and
all pinned vertices of out-degree $0$), the strongly connected
decomposition (with the pinned vertices identified) coincides with
the \blocktriangular decomposition of the pinned rigidity matrix
with a maximal number of diagonal blocks, for some linear order of the blocks
extending the partial order of the directed graph decomposition. \label{thm:matrixdecom}
\end{theorem}

\proof  Given a pinned $d$-isostatic graph with a pinned
$d$-directed orientation, we apply the techniques of \S2 to the
pinned graph with the pinned vertices identified to obtain a
strongly connected decomposition. If this has more than one
component, then focus some  bottom strongly connected component $A$, extended
with its edges to the ground.  With a permutation of the
rows and a permutation of the columns, we place these vertices and
edges at the top left of the matrix.  The rest of the rows and columns for a second block $B$ (that may not be strongly connected).  (For notation in the matrices below, we assume with
$i,j\in A$ and $\ell,m \in B$ and $k$ is a pinned vertex, where
$B$ is the  component.) This gives the form:

\noindent\makebox[\textwidth]{%
\begin{tabularx}  {\textwidth}{X}
\[ \begin{array}{c|ccccc|ccccc|}
 &  \ldots  & i&  \ldots& j & \ldots&\ldots & \ell& \ldots& m & \ldots\\
\hline  &  \ddots &  &  \vdots& & \Ddots&   \ddots & &   \vdots  &  & \Ddots \\
\{i,j\} &  \ldots 0&p_{i}-p_{j}& 0 \dots0&p_{j}-p_{i}& \ldots& 0 \ldots  & 0&   \ldots  &0 & \ldots  0\\
\{i,k\} &  \ldots 0&p_{i}-p_{k}&  0\dots\  \ &0& \dots& 0 \ldots & 0&  \ldots &0 &  \ldots 0\\
 &  \Ddots &  &  \vdots&  & \ddots& \Ddots  &  &  \vdots &  & \ddots\\
\hline
 &  \ddots &  &  \vdots&  & \Ddots& \ddots & & \vdots  & & \Ddots\\
\{i,\ell\} &  \ldots 0& p_{i}-p_{\ell}&  0\ldots&0& \ldots& \ldots 0&p_{\ell}-p_{i}&  0 \ldots 0 &0&0 \ldots\\
\{\ell,m\} &   \ldots 0& 0&\ \  \ldots& 0& \ldots& \dots 0 &p_{\ell}-p_{m}&0 \ldots 0&p_{m}-p_{\ell}& 0\dots\\
 &  \Ddots &  &  \vdots&  & \ddots& \Ddots & &  \vdots& &\ddots\\
\hline
\end{array}  \]
\end{tabularx}}
or simply:
\begin{displaymath}
\mathbf{R}(G) =
 \begin{tabularx}{.1\textwidth}{X}
\[ \begin{array}{|c|c|}
\hline
  \mathbf{R}(A) &0\\
\hline
  X &\mathbf{R}(B)\\
\hline
\end{array}  \]
\end{tabularx}  \phantom{xxxxxxxxxxxxxxx}
\end{displaymath}
We have an initial block decomposition.  Repeating this for each
of the blocks up the acyclic strongly connected decomposition, we
find a diagonal matrix block for each component of the
decomposition, and the matrix has the desired \blocktriangular
form (with $B_{2}, ... B_{n}$ as the remaining blocks):
\begin{displaymath}
\mathbf{R}(G) =
 \begin{tabularx}{.1\textwidth}{X}
\[
\begin{array}{|c|c|c|c|c|}
\hline
  \mathbf{R}(A) &0& 0&\ldots & 0\\
\hline
   X_{12} &\mathbf{R}(B_{2})&0& \ldots & 0\\
\hline
   \vdots &\vdots&\ddots &\ddots & \vdots\\
\hline
   X_{1n} &X_{2n}& X_{3n}&\ldots &\mathbf{R}(B_{n})\\
\hline
\end{array}
\]
\end{tabularx}
\phantom{xxxxxxxxxxxxxxxxxxxxxxxxxxxxxx}
\end{displaymath}

Note that the acyclic strongly connected decomposition graph only
gives a partial order.  The matrix  \blocktriangular decomposition
imposes a linear order which extends this partial order.  The linear order will
not be unique, as it is possible that some blocks may be
incomparable in the partial order - which will show up here as
some off diagonal blocks with $X_{ij}=0$.  (See Figure~\ref{Block} where component $(GHIJ)$ could be moved anywhere to the right in the linear order.)

Conversely, if we have such a block-triangular decomposition, then
the corresponding Laplace decomposition of the determinant will
produce a pinned $d$-directed orientation in which all edges are
either oriented within the block, or directed towards vertices in
blocks to the left.  In short, we have a decomposition with an
acyclic decomposition graph extracted from the linear order of the
blocks. 

If there is a further decomposition from the $d$-directed graph, then
this will show up as a further block-triangular decomposition with
more components, and vice versa. \eop


Recall that any two pinned $d$-directed orientations are
equivalent, so they give the same decomposition in the algorithm of
\S2. Therefore in the \blocktriangular decomposition of the pinned
rigidity matrix, the linear order extends this partial order.  In
particular, each other non-zero term in the Laplace block
expansion gives an equivalent pinned $d$-directed orientation, and
an equivalent \blocktriangular decomposition of the pinned
rigidity matrix.  Moreover, it is easy to check that any pinned
$d$-directed orientation of a graph can be used to create a non-zero term in
the Laplace decomposition of the determinant. There is a bijection
between these $d$-directed orientations and the non-zero terms of the Laplace
expansion.

\noindent {\bf Remark.}  We note that CAS programs such as Maple,
Mathematica and SAGE have `strongly connected decompositions of
matrices' into lower triangular form \cite{Monagan}.   However,
these are not the same - only analogous.  The `directed graph'
used in their algorithm comes from treating the rows and columns of
the square matrix as indexed by the same vertex set, and using any
non-zero entry $a_{i,j}$ as a directed edge $(i,j)$.  The strongly
connected decomposition then corresponds to a simultaneous
permutation of the rows and columns, so that we go from $A$ to the
similar matrix $[D]=P[A]P^{-1}$.   In our process, we have
identified a different graph with fewer vertices ($|V|$ instead of
$d|V|$) and different edges.  Our shift allows  distinct
permutations of the rows and columns,  generating a more general
$[D']=P[A]Q$.  With more possible permutations, we get all the
decompositions that flow from the standard CAS code - and more. In
particular, a simple permutation of the rows of the pinned
rigidity matrix will change a matrix which Maple declares does not
decompose, into one that Maple notices does decompose!   \eop


\subsection{$d$-Assur graphs and minimal components}\label{subsec:2matrix}

We now focus on the minimal components in these decompositions of a  pinned $d$-isostatic graph
$\widetilde{G} = (I,P;E)$.  We will verify that these are minimal
pinned  $d$-isostatic components of $\widetilde{G}$, i.e. the {\it
\dassur graphs}. We also verify some key properties of these graphs.

\begin{theorem}[\dassur Graphs] For a pinned  $d$-isostatic graph $\widetilde{G}$ the following are equivalent:
\begin{enumerate}
  \item the graph contains no proper pinned $d$-isostatic subgraphs;
  \item the graph is indecomposible for some (any)  $d$-directed orientation;
  \item the pinned rigidity matrix has no proper block triangular decomposition.
\end{enumerate}\label{thm:2blockdecomp}
\end{theorem}

\proof  Theorem~\ref{thm:matrixdecom} shows the equivalence of (2) and (3).   It remains to show that (1) is equivalent to these.

Assume the graph is not minimal and there is a proper pinned isostatic subgraph $\widehat{G}$.  We will show the pinned rigidity matrix decomposes.   If we permute the inner vertices and all the edges of $\widehat{G}$ to the upper left corner of the pinned rigidity matrix, the rest of these rows from $\widehat{G}$ are $0$ and the remaining columns and rows form a second block.  This gives a block triangular decomposition of the pinned rigidity matrix.   The contrapositive says that if the pinned rigidity matrix does not have a proper block triangular decomposition, then the pinned isostatic graph is minimal.

Conversely, if the pinned rigidity matrix for $\widetilde{G}$   has a proper block triangular decomposition, then it is clear the upper left block represents the inner vertices and the edges (including the ground edges) of a proper pinned isostatic subgraph  $\widehat{G}$.  The contrapositive confirms that if the graph  $\widehat{G}$ is minimal then the pinned rigidity matrix does not decompose.
 \eop

 Any graph $\widetilde{G}$ which satisfies one of these three equivalent properties will be called \dassurp.
 We can use the concept of a minimal pinned isostatic graph (\dassur graph)  to build up a third route to the decomposition of a larger pinned $d$-isostatic graph into \dassur graphs.  We follow the same track used for \tassur graphs in \cite{SSWI}.

\begin{figure}[htb]
\centering
 \subfigure[] {\includegraphics [width=.4\textwidth]{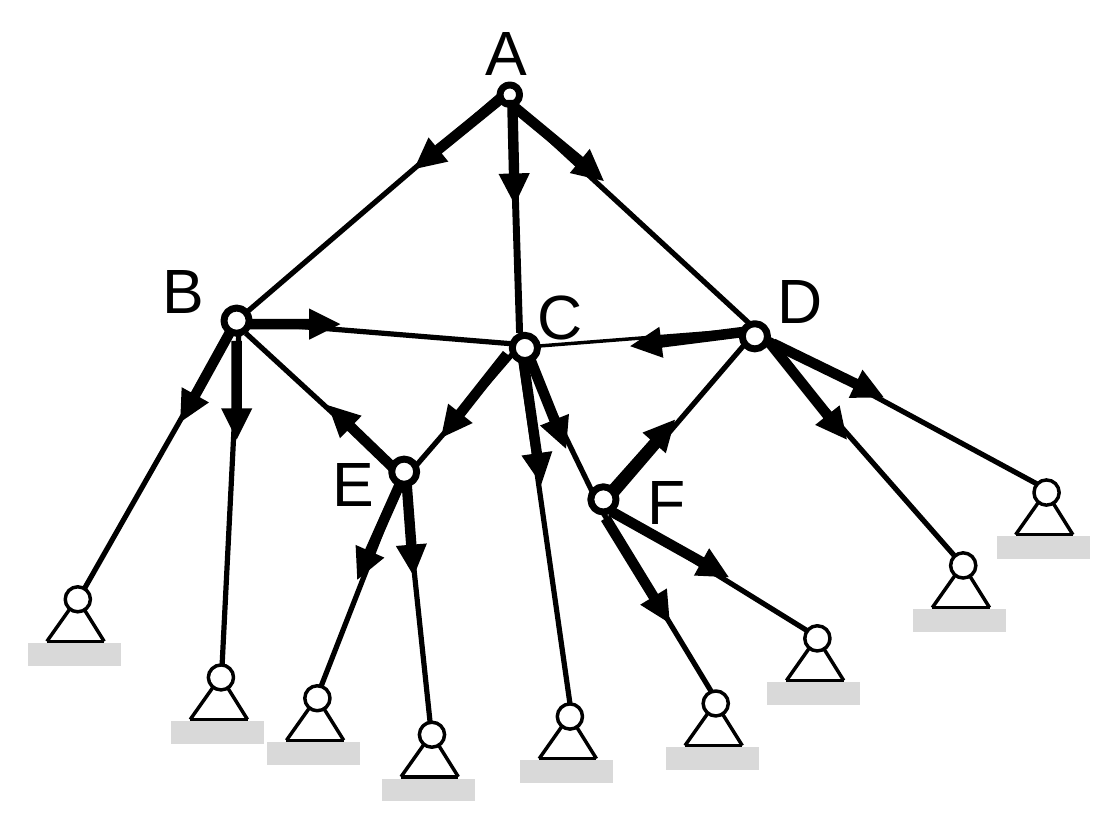}} \quad
   \subfigure[] { \includegraphics[width=.3\textwidth]{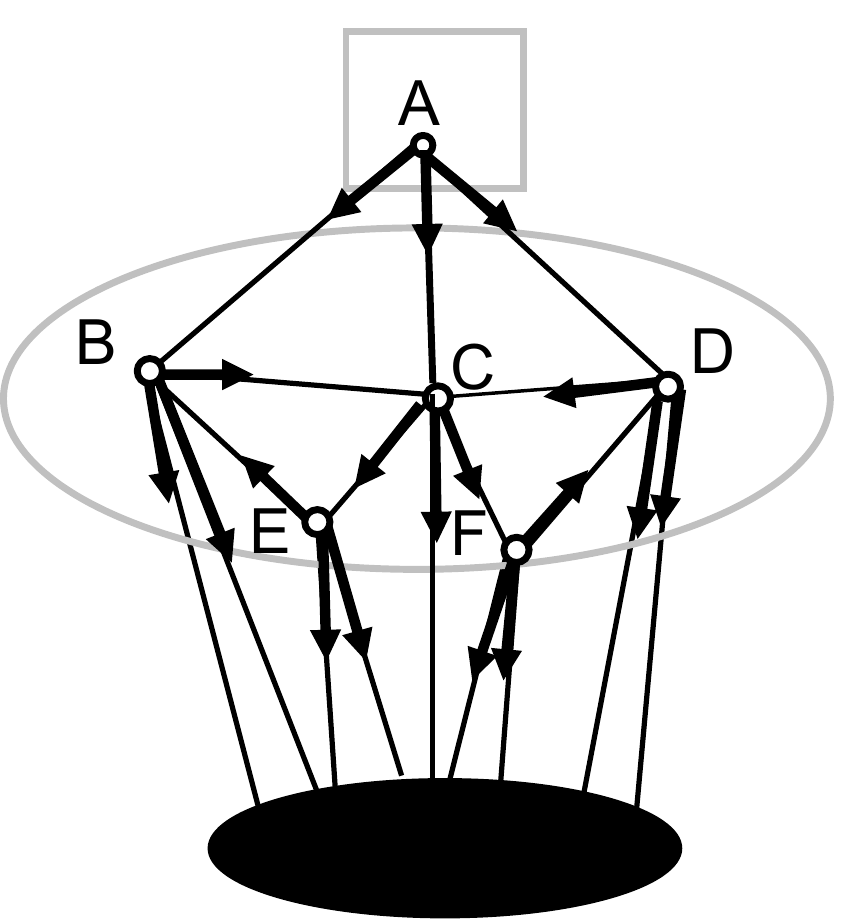}} \quad
      \subfigure[] { \includegraphics[width=.25\textwidth]{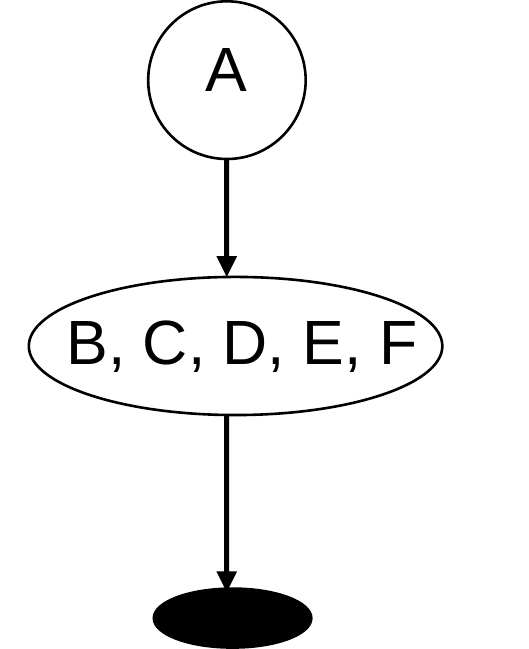}}\quad
          \subfigure[] { \includegraphics[width=.45\textwidth]{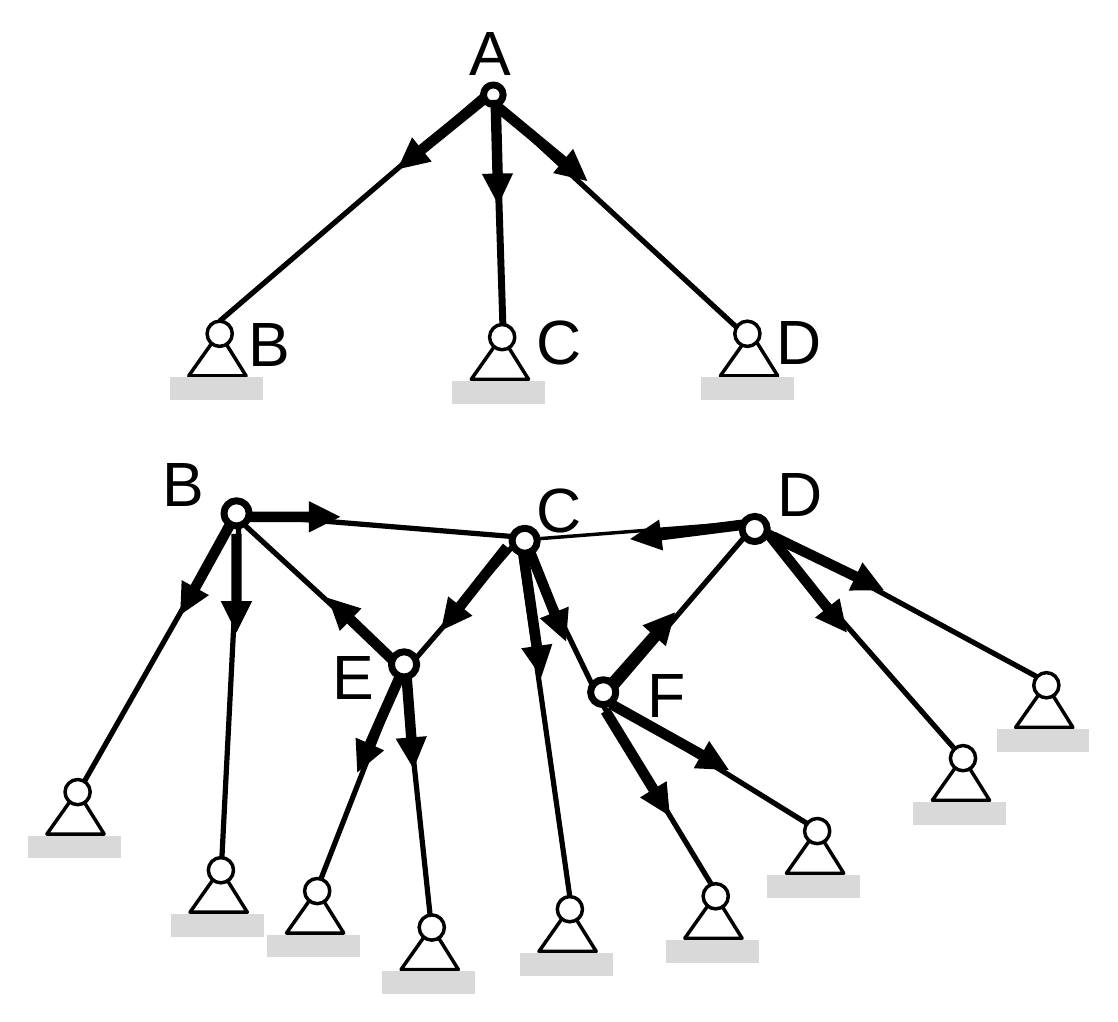}}
\caption{A pinned $3$-isostatic framework (a) has a unique
decomposition into 3-Assur graphs through condensing the pins to
ground (rigid region) (b) and creating a partial order (c) which
also corresponds to a scheme of Assur graphs (d).\
\label{fig:3dassurgraph}}
\end{figure}

 We just recall the outline of the process, since the decomposition will match what we already have from Theorem~\ref{thm:matrixdecom}.  (See Figure~\ref{fig:3dassurgraph}.)
The ground is the bottom layer.  Find a minimal $d$-isostatic pinned subgraph.  This will be above the ground component.  Combine that into the ground, and seek another minimal $d$-isostatic pinned subgraph.  This will be  a component above all components to which its ground edges attach.     Repeat until all vertices are in some component.  This is the {\it \dassur decomposition} of the pinned $d$-isostatic graph $\widetilde{G}$.

\subsection{Summary Assur Decomposition} \label{AssurDecomp}
We now pull these connections into a summary theorem about the three equivalent ways to decompose a pinned $d$-isostatic graph.

\begin{theorem}[\dassur Decomposition]  Given a pinned $d$-isostatic graph $\widetilde{G}$, there is a $d$-directed orientation of $\widetilde{G}$. With any such $d$-directed orientation, the following
decompositions are equivalent:
\begin{enumerate}

\item  the \dassur decomposition of $\widetilde{G}$;

\item   the strongly connected decomposition into extended
components associated with the $d$-directed orientation (with all
pins identified);

\item  the  \blocktriangular  decomposition of the pinned rigidity
matrix into a maximal number of components for some linear order
extending the partial order of (i) or equivalently (ii).

\end{enumerate}
\end{theorem}

\proof  Theorem~\ref{thm:matrixdecom} give the equivalence of (2) and (3).  The equivalence of (1) and (2)  follows from Theorem~\ref{thm:2blockdecomp}.  The  construction process for the $d$-Assur decomposition above guarantees that for two components $A,B$ $A$ is above $B$ in the  $d$-Assur  partial order if and only if $A$ is above $B$ in the  $d$-directed  partial order.
\eop

\subsection{Pinned $d$-counts and insufficiency of $d$-directions}\label{subsec:counts}
While having a $d$-directed orientation is a necessary condition
for being $d$-isostatic, it is not sufficient.  These orientations are also connected to some basic necessary counting conditions - conditions that have been proven sufficient in the plane, but fail to be sufficient in $3$-space.  Since there have been efforts to transform this necessary condition into a sufficient condition, we give a few key observations and examples to clarify the connections.

\begin{figure}[htb]
\centering
 \subfigure[] {\includegraphics [width=.55\textwidth]{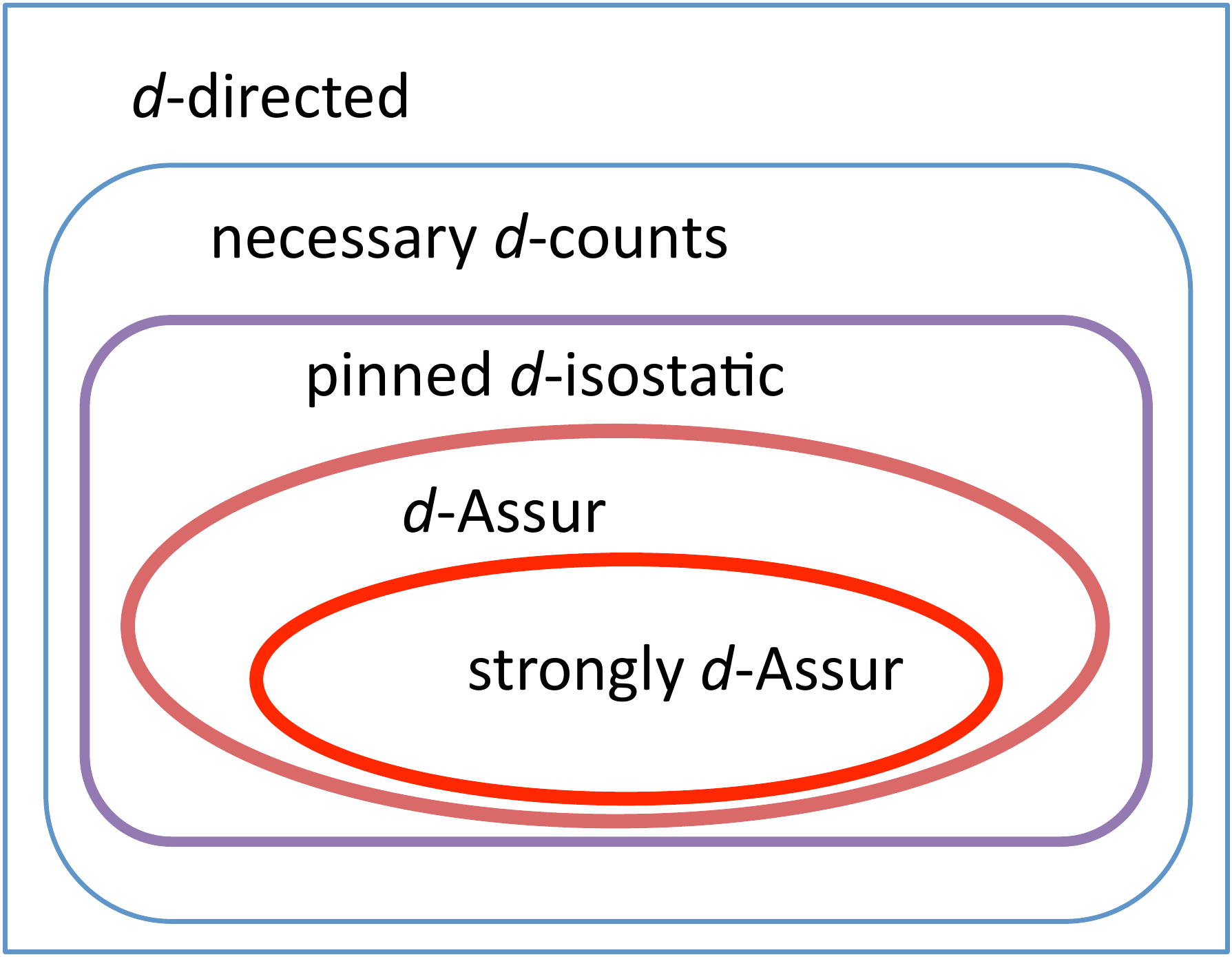} }\quad
 \subfigure[] {\includegraphics [width=.40\textwidth]{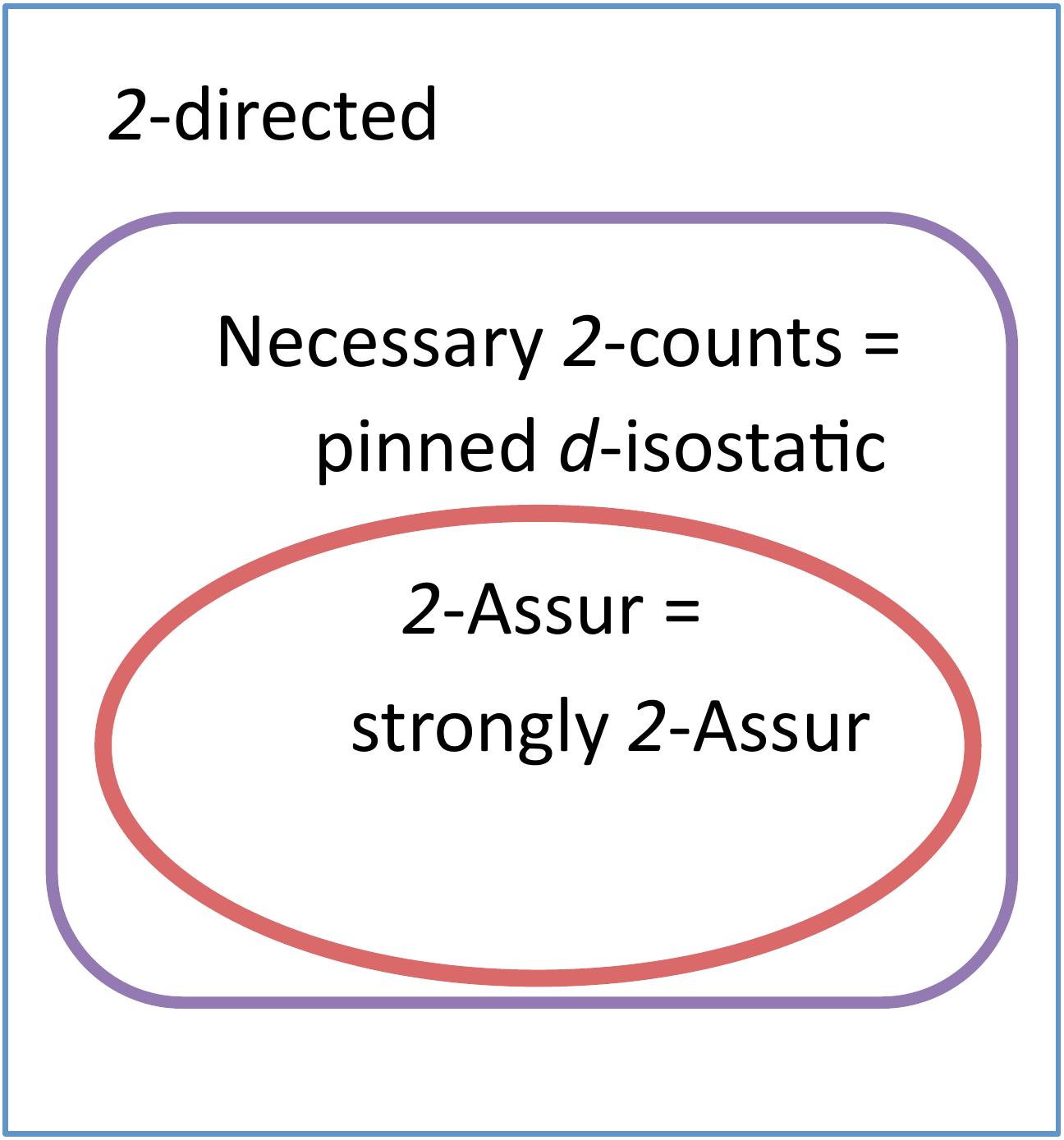} }
\caption{(a) shows the general nesting of conditions for dimension $d$, including \sdassur graphs, while (b) shows the simplified nesting in the plane.\label{fig:assurweaklyassur}}
\end{figure}

\begin{theorem}[Necessary Pinned $d$-Counts]\label{3PinThm}
Given a pinned $d$-isostatic graph $\widetilde{G} = (I,P;E)$, the
following properties hold

\begin{enumerate}
\item[(1)] $|E|=d|I|$.

\item[(2)] for all subgraphs $\widetilde{G}' = (I',P';E')$
\begin{itemize}
\item[(i)]         $|E'|\leq d|I'|$ if $|P'| \geq d$,

\item[(ii)]          $|E'|\leq d|I'|-{d + 1 - k\choose 2}$ if $|P'
|= k$,  $0 \leq k< d$,
\end{itemize}\end{enumerate}

\end{theorem}

\proof  The count (1) is just the condition for the pinned
$d$-rigidity matrix to be square.  For the subgraphs, the count
2(i) is again the subgraph count on the submatrix of the columns
of $I'$, whose failure would guarantee a row dependence.

In 2(ii) ${d + 1 - k\choose 2}$ is the count of remaining
rotational degrees of freedom with k points pinned. For instance,
with $d = 3$ and $k = 2$ ($|P' |= 2$), there must be a ${2 \choose
2} =  1$ remaining degree of freedom, corresponding to a rotation
of the entire framework about a line passing through these two
points. Note that pinning one point eliminates all d translational
degrees of freedom. As we continue to pin additional vertices we
are eliminating more rotational degrees of freedom. With $k$
points fixed, $k - 1$ dimensional space is fixed. The rotational
axis in any dimension $d$ is of dimension $d -2$ (i.e. point in
plane, line in $3D$, plane in $4D$, etc). The possible rotations
occur in the orthogonal complement to the rotational axis, which
is a $2$-dimensional space. With $k$ points fixed from the
remaining $d - (k-1)$ space we choose all such possible
$2$-dimensional orthogonal complements which gives us the
remaining ${d - (k-1)\choose 2}$ $=$ ${d + 1 - k\choose 2}$ space
of rotations (DOF). For 2 (iii), there is a space of ${d \choose
2}$ rotations fixing the single vertex in dimension d, and for
2(iv) there is a ${d+1 \choose 2}$-space of trivial motions
(combination of $d$ translations and ${d \choose 2}$ rotations),
so the bounds on the count of edges is the maximum number of
independent rows. \eop

\medskip

Note that the formula for the subgraph count in 2(ii) of
Theorem~\ref{3PinThm} is also the correct form for all $k$ $\geq$ $0$.

In dimension 2, the necessary pinned 2-counts become ({\em Pinned
Plane Framework Conditions} for $\widetilde{G}=(I,P;E)$)
\cite{SSWI}:

\begin{enumerate}
\item $|E|=2|I|$ and \item for all subgraphs $\widetilde{G}'=(I',
P'; E')$ the following conditions hold:
\begin{itemize}
\item[(i)]          $|E'|\leq 2|I'|$ if $|P'| \geq 2$,

\item[(ii)]            $|E'|\leq 2|I'| -1$ if $|P'|=1$ ,  and

\item[(iii)]            $|E'|\leq 2|I'|-3$ if $P' = \emptyset$.
\end{itemize}
\end{enumerate}

For completeness, the Necessary pinned 3-counts are:

\begin{enumerate}
\item[(1)] $|E|=3|I|$.

\item[(2)] for all subgraphs $\widetilde{G}' = (I',P';E')$
\begin{itemize}
\item[(i)]         $|E'|\leq 3|I'|$ if $|P'| \geq 3$,

\item[(ii)]          $|E'|\leq 3|I'| -1$ if $|P'|=2$ ,  and

\item[(iii)]          $|E'|\leq 3|I'|-3$ if $|P' |= 1$.

\item[(iv)]          $|E'|\leq 3|I'|-6$ if $P' = \emptyset$.
\end{itemize}\end{enumerate}

A useful observation is that all pinned graphs in dimension $d$
that satisfy the necessary pinned $d$-counts are $d$-directed:

\begin{theorem}\label{pinnedcountsdirected}
If a pinned graph $\widetilde{G} = (I,P;E)$ satisfies the
Necessary pinned $d$-counts for some $d$, then it has a
$d$-directed orientation with all inner vertices of out-degree $d$
and inner vertices of out-degree $0$.
\end{theorem}

One proof is to use the pebble game algorithm operations,
essentially playing a $d|V|$ pebble game on both inner edges and
ground edges, giving us the desired $d$-directed orientation of
the graph \cite{mechanicalpebblegame}.

\begin{figure}[htb]
\centering
   \subfigure[] {\includegraphics[width=.38\textwidth]{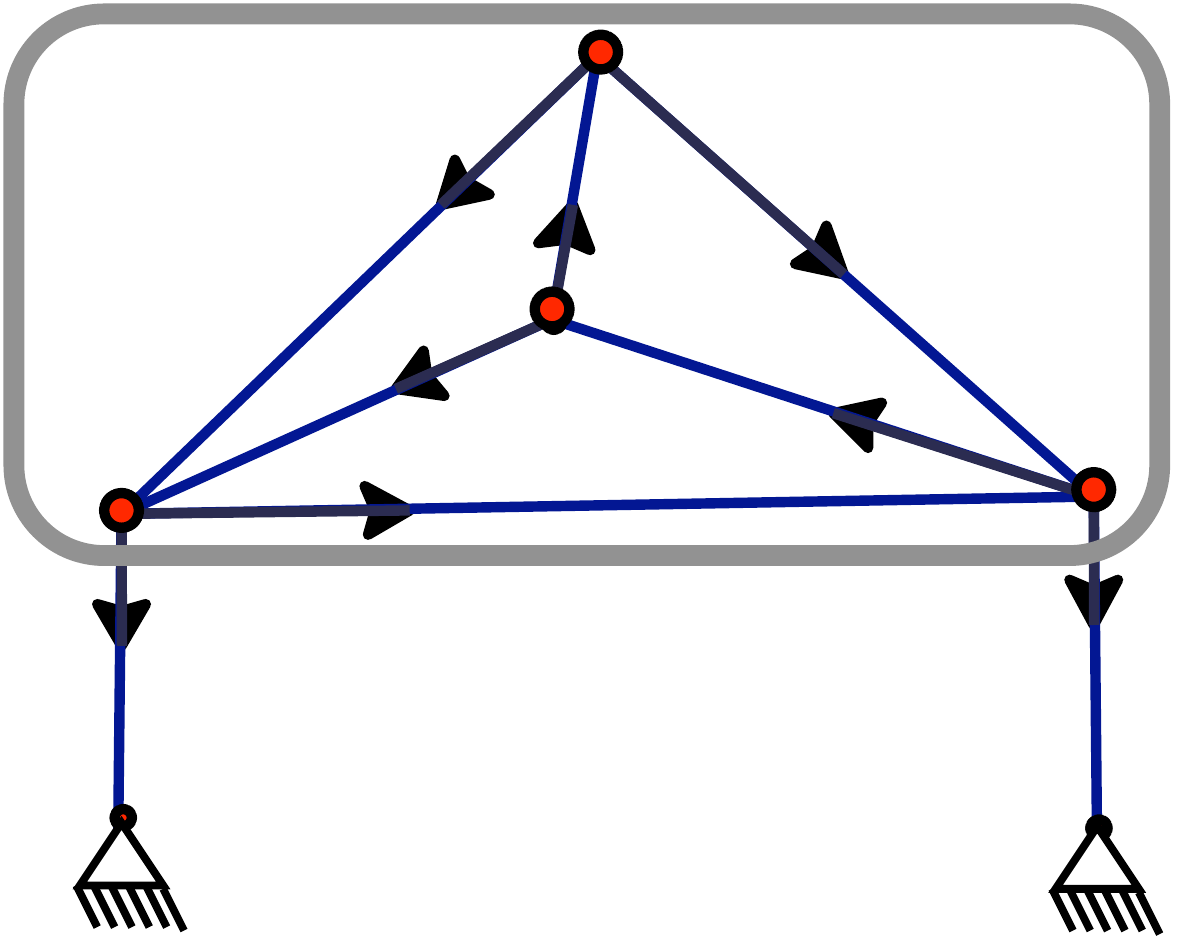}\qquad\quad
   \subfigure[] {\includegraphics [width=.30\textwidth]{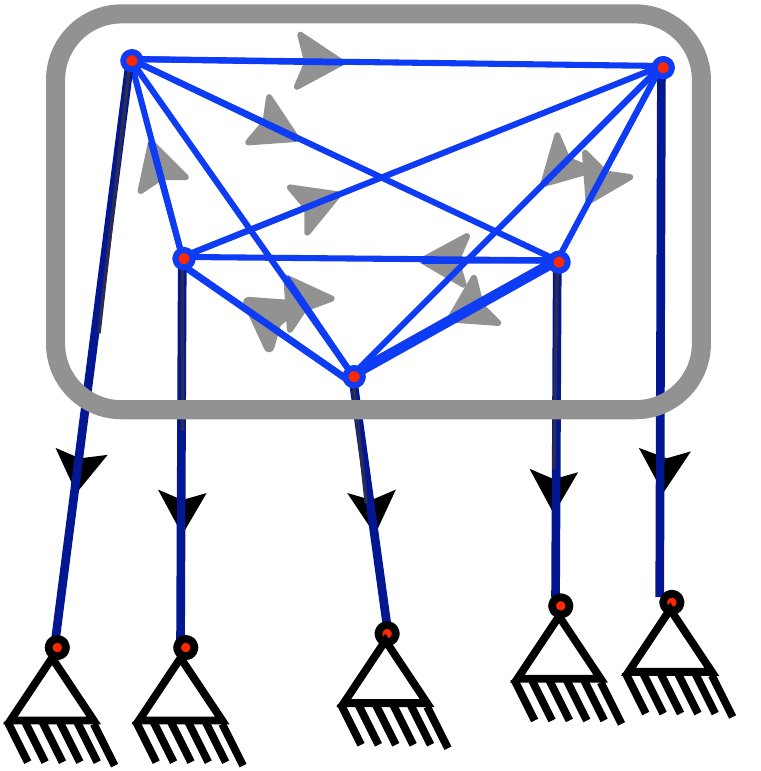}}}
\caption{A pinned graph which is $2$-directed, indecomposible but
not 2-isostatic (the subgraph in the grey box is overcounted - it
does not satisfy the necessary pinned $2$-counts). (a). A pinned
graph which is 3-directed, indecomposible but not 3-isostatic (the
subgraph in grey box does not satisfy the necessary 3-counts) (b).
\label{2bad}}
\end{figure}

Having stated the necessary pinned $d$-counts, we can now
illustrate that $d$-directed orientation is not a sufficient
condition for being $d$-isostatic. In Figure~\ref{2bad} (a) we
have an example of a graph which is $2$-directed but not
$2$-isostatic. It is not $2$-isostatic as it fails the necessary
pinned $2$-counts. Such a graph will have a non-zero term in the
Laplace expansion, and a block decomposition. However the
determinant is zero, overall, with cancellation among non-zero
terms. The example in Figure~\ref{2bad} (b) shows a similar
example in 3-space, it has a proper $3$-directed orientation,
although it is not pinned $3$-isostatic. These two examples
confirm that having a 2 and 3-directed orientation and being
indecomposible does not capture the very basic subgraph counting
conditions of Theorem~\ref{3PinThm}.

\begin{figure}[htb]
\centering
  \subfigure[] {\includegraphics [width=.3\textwidth]{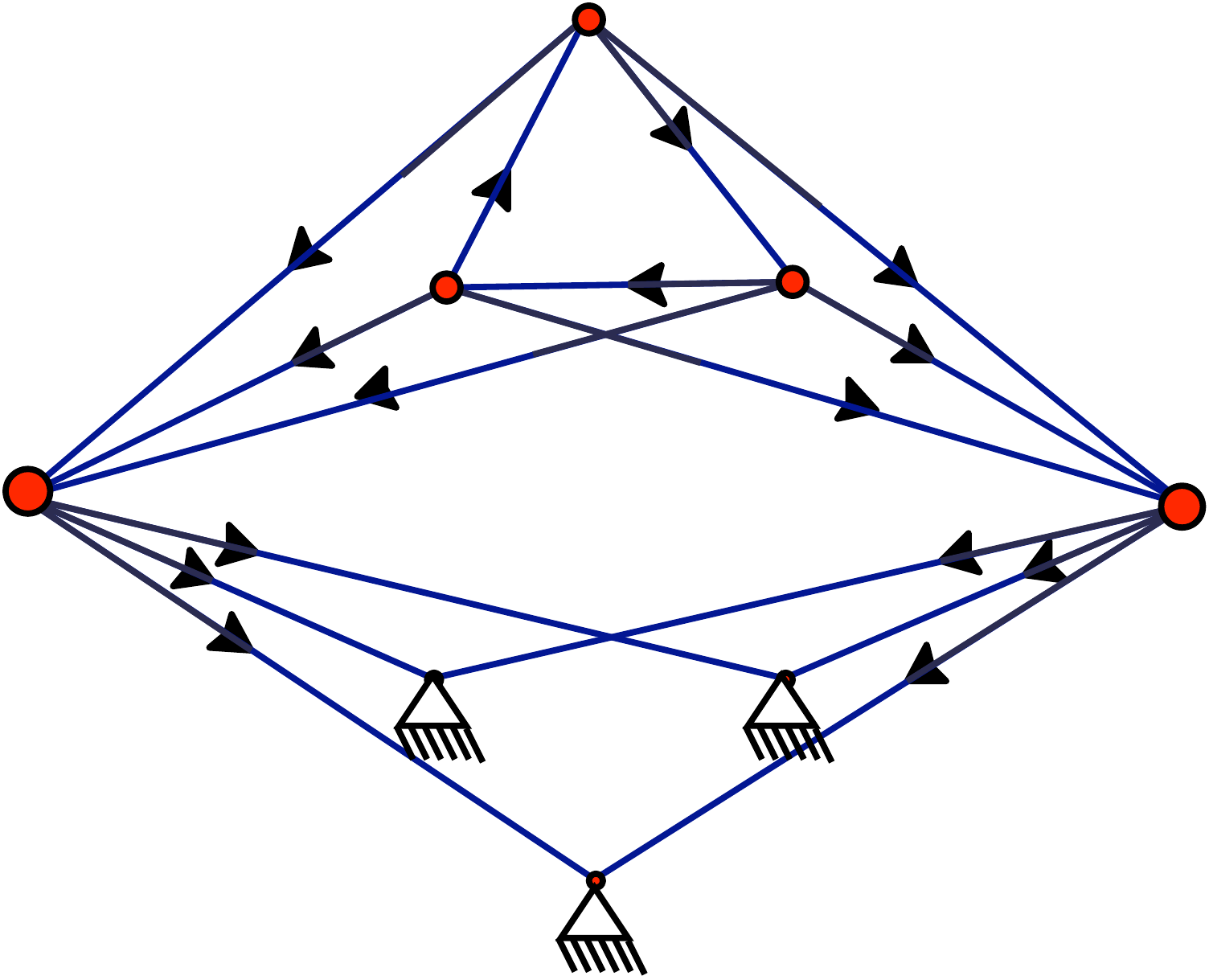}} \quad
    \subfigure[] {\includegraphics [width=.35\textwidth]{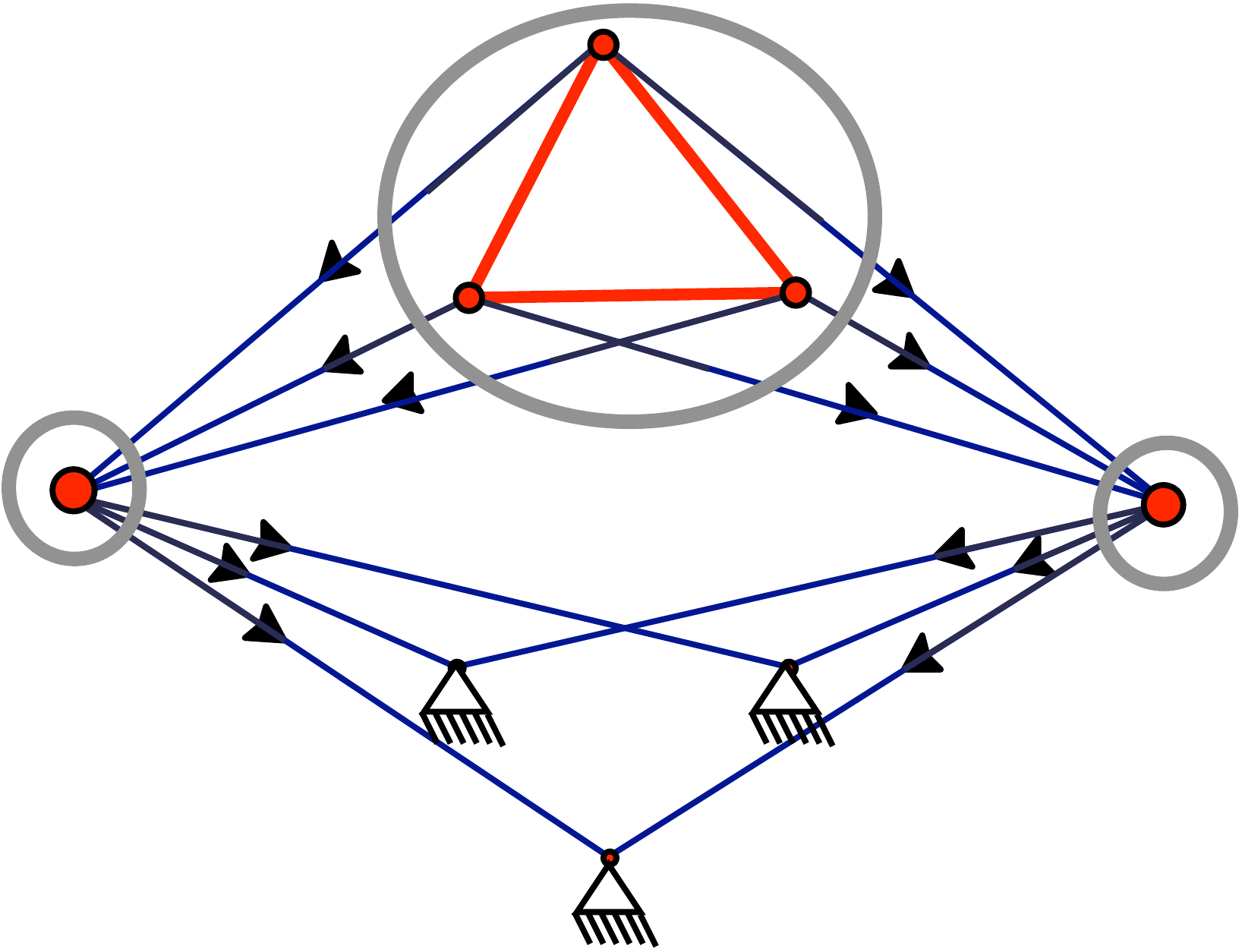}} \quad
 \subfigure[] { \includegraphics[width=.25\textwidth]{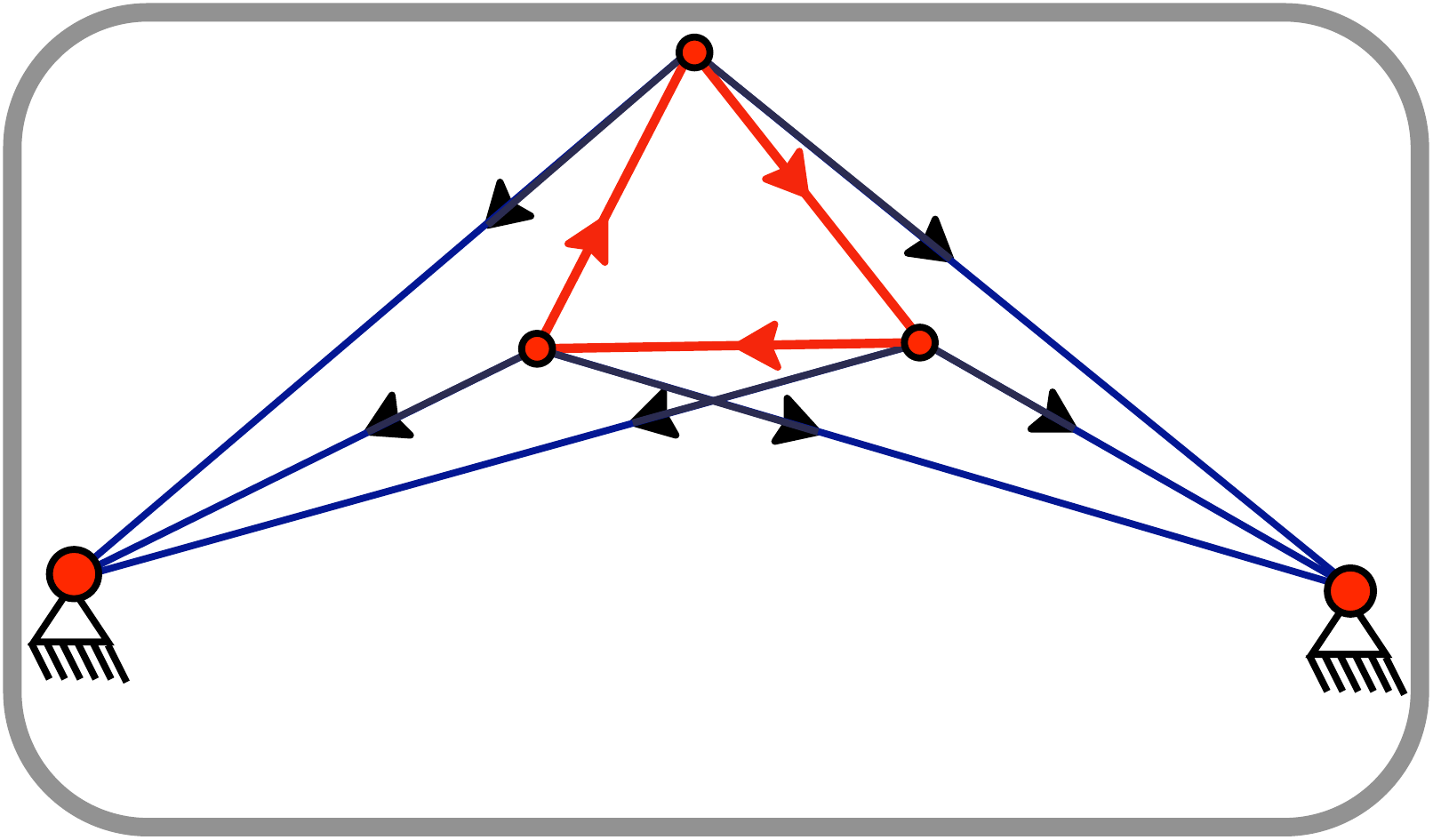}}
 \caption{A pinned graph graph which is $3$-directed (a), has a strongly connected decomposition (b)
 but the top extended component (c)  is not 3-isostatic.\label{banana2}}
\end{figure}

In Figure~\ref{banana2} we have another example which has a
3-directed orientation (a), although it is not pinned 3-isostatic
(see remark below for further discussion on this special example).
This graph is decomposible (b), but what is interesting is that
one of the components (c) now visibly fails the necessary pinned
3-count though is still $3$-directed.

\noindent {\bf Remark.} By translating the general results for
generically isostatic graphs in dimension 1 and 2, the pinned
counts in dimension 1 and 2 are also sufficient for a graph to be
pinned 1 and 2-isostatic, respectively. In dimension 2, the
necessity and sufficiency of the counts is captured in the Pinned
Laman Theorem \cite{SSWI}, which states that $\widetilde{G}$
satisfies {\em Pinned Plane Framework Conditions} if and only if
it there exists a pinned 2-isostatic realization of
$\widetilde{G}$ in the plane.

In dimensions $d > 2$ however, these counting conditions are not
sufficient for $d$-isostatic graphs. There are classic examples
which show that these counts are not sufficient in $3$-space.
Figure~\ref{banana2} is one such classic example (analogue to a
well known `double banana' example in 3D unpinned bar and joint
frameworks, also known in mechanical engineering community as
floating frameworks). \eop

A $d$-directed graph has an immediate directed graph decomposition.  Whether the graph is pinned $d$-isostatic will depend on whether all of the components are pinned $d$-isostatic.

\begin{figure}[htb]
\centering
 \subfigure[] {\includegraphics [width=.55\textwidth]{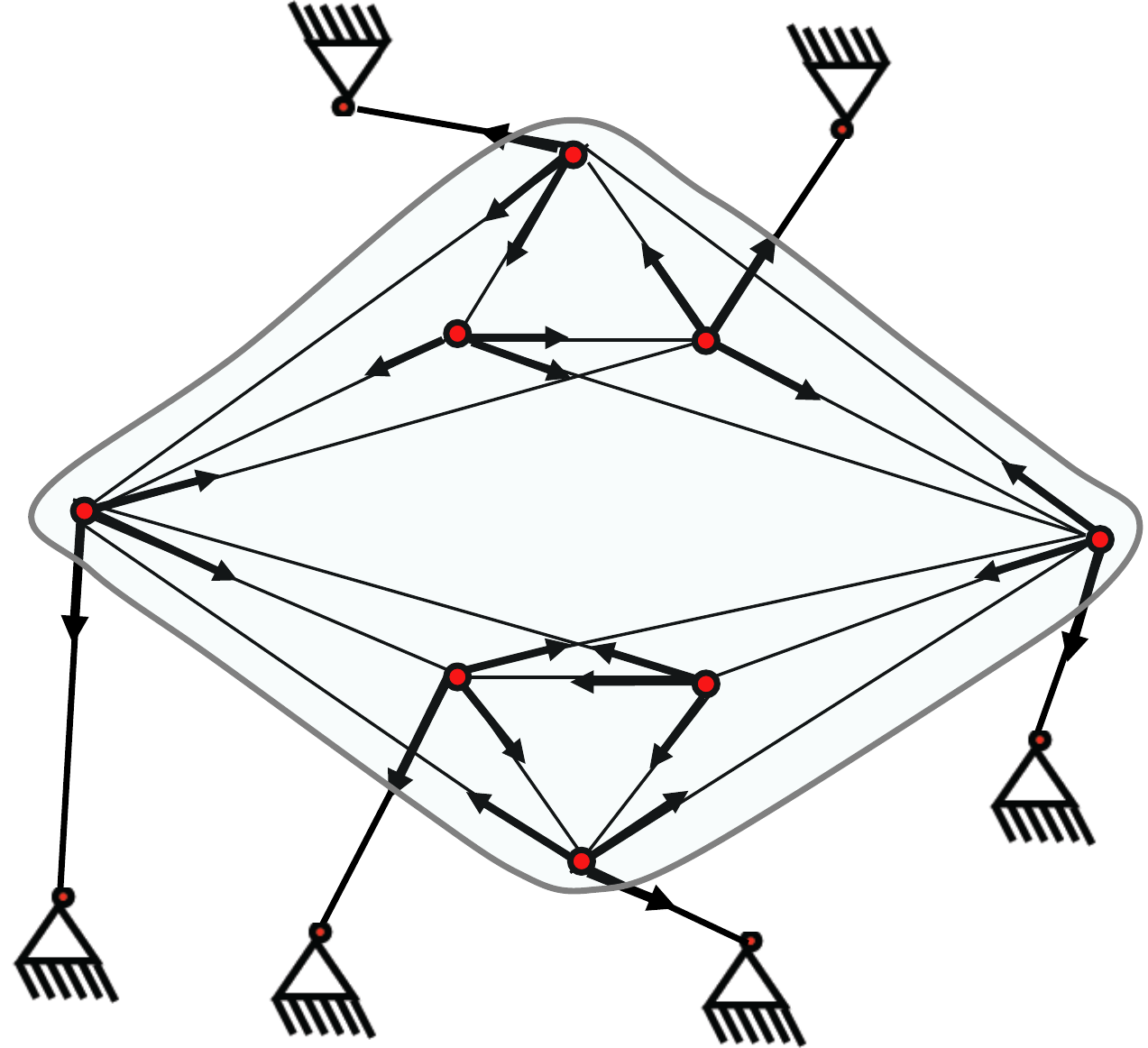}}
 \caption{Pinned graph which is $3$-directed, indecomposible, satisfies the necessary pinned $3$-count, but it is not 3-isostatic.\label{banana3}}
\end{figure}

While such a directed graph decomposition detects some failures, via failed
necessary counts on subgraphs, detecting from a graph if a
$d$-pinned framework is $d$-isostatic for $d$ $>$ $2$ is generally
difficult. Alternatively, we would have to resort to analysis of
the pinned rigidity matrices. Figure~\ref{banana3} shows a
$3$-directed indecomposable graph which satisfies all the subgraph
necessary pinned $3$-counts of Theorem~\ref{3PinThm} but is still
not pinned $3$-isostatic. Even combined with the \dassur
decompositions, we do not have necessary and sufficient counting
conditions for a graph to be \dassur when $d$ $>$ $2$.
There are simple algorithms to detect the failure of the type in
Figure ~\ref{2bad}, but no known polynomial algorithms for the
failures of type Figure~\ref{banana3} \cite{LeeStreinu,
mechanicalpebblegame}.

In 2D it is easier to check whether any pinned graph is
$2$-Assur (Figure~\ref{2bad}(b)). We need to test both (i) the complete set of Pinned
Plane Framework Conditions and (ii) the indecomposibility.  While
the decomposition algorithm is linear in $|E|$, the process of
extracting the $2$-directed orientation from the pinned rigidity
matrix is exponential. In \cite{mechanicalpebblegame} we have
presented a polynomial time ($O(|I|^3)$) algorithm to test both
(i) and (ii). This algorithm is based on the strongly connected
decompositions of the graph and the pinned version of the pebble
game algorithm.

However, if we already know that the graph is pinned
$d$-isostatic, regardless of the dimension $d$, we can use the
pinned pebble game algorithm to directly check whether the graph is
\dassur or if it can be further decomposed to individual
\dassur components \cite{mechanicalpebblegame}.

\section{Tracing Motions in Linkages through Assur Decompositions}\label{sec:flex}
In the previous papers for plane Assur graphs \cite{SSWI,SSWII},
some additional properties were explored.  In the context of
linkages, a central property is how an associated framework
responds when one of the links (edges) is
replaced by a `driver' so that this distance between the end
points is controlled by a piston, or we control an angle at a hinge, as in a robot arm.

\begin{figure}[htb]
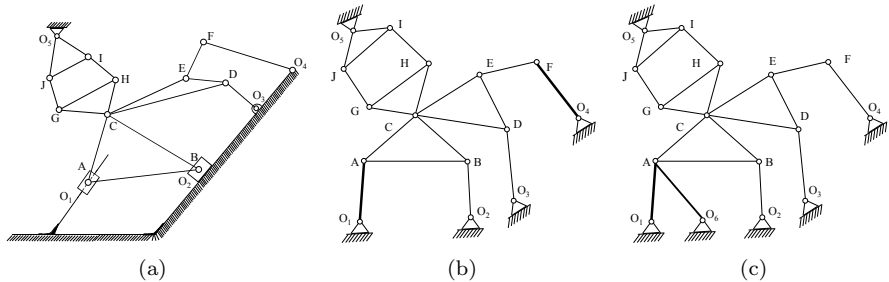

\centering
  \subfigure[] {\includegraphics [width=.33\textwidth]{Linkage.pdf}} \
  \subfigure[] {\includegraphics [width=.31 \textwidth]{StructuralScheme.pdf}}
  \subfigure[] {\includegraphics [width=.31 \textwidth]{AdaptedScheme1.pdf}}
 \caption{A plane linkage (a) is translated to a framework scheme with one degree of freedom (b) which is then made pinned $2$-isostatic by adding a bar to block a chosen driver (c) .\label{fig:edgeremoval}}
\end{figure}

We first  recall how some types of motions of linkages
were translated into pinned frameworks.  Then we consider the
infinitesimal motions of inner vertices which arise when one edge
is removed.  This is followed by some additional variations in how
`drivers' are inserted into pinned isostatic frameworks. We will
word the theorems whenever possible for \dassur graphs,
though several key results only hold for $2$-Assur graphs.

In \S4.2 we will distinguish the special class of {\it \sdassur
graphs} by the feature that all inner vertices go into motion for
every choice of inserted driver. This is the direct extension of a
stronger key property of $2$-Assur graphs \cite{SSWI}.
This is not a matter of a deeper decomposition - but a difference in the characteristics of the underlying components we are analyzing (or synthesizing).

\subsection{From linkages to structural schemes}
In the introduction, we presented a figure of a linkage, complete
with some slider joints, side by side with the {\em structural
scheme} - corresponding flexible pinned bar and joint framework.
The results in \S3 applied to pinned isostatic frameworks, and
those in \S2 applied whenever we generated a directed graph.

\begin{figure}[htb]
\centering
 \subfigure[] {\includegraphics [width=.34\textwidth]{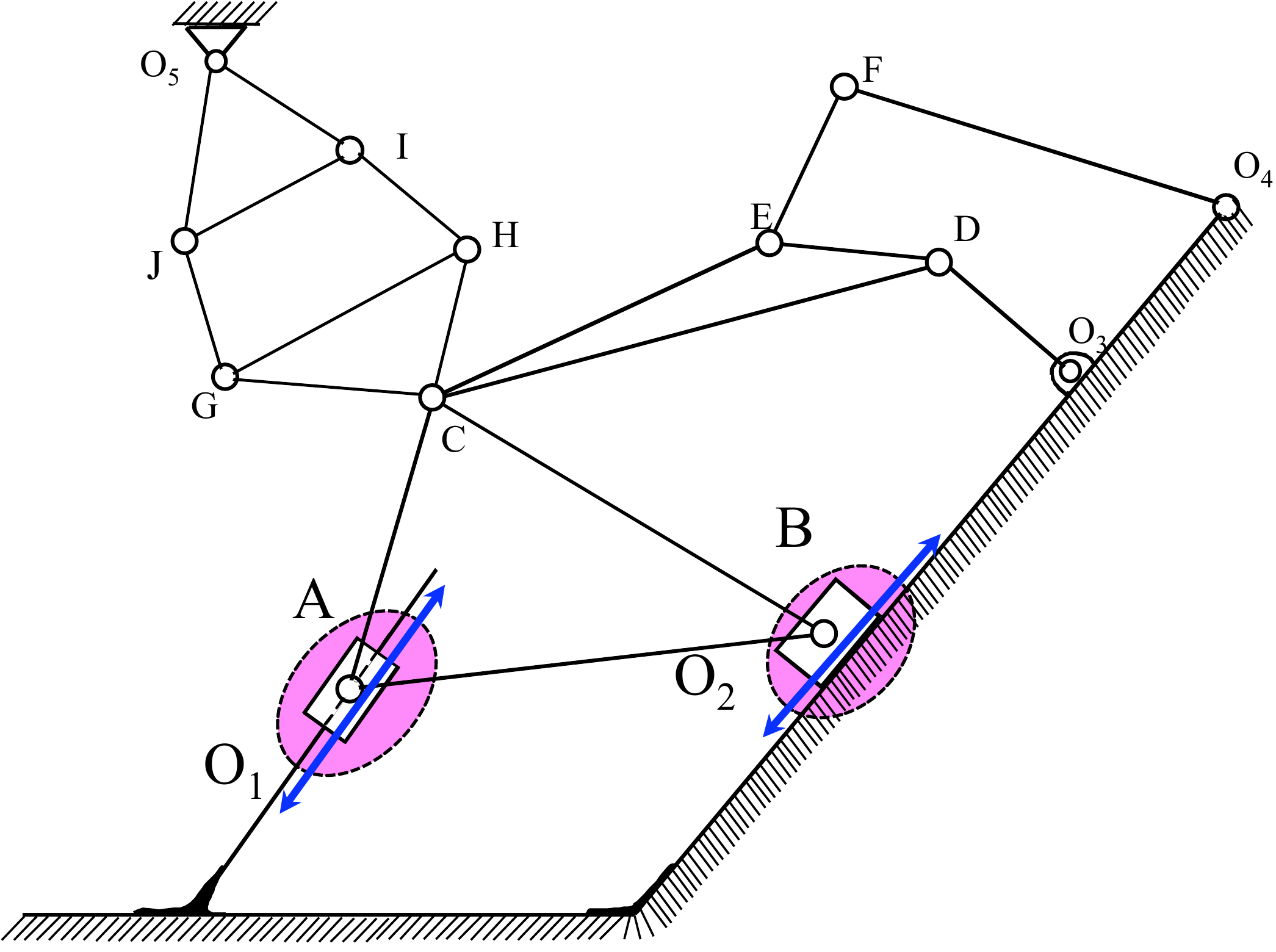}}\!\!
  \subfigure[] {\includegraphics [width=.34\textwidth]{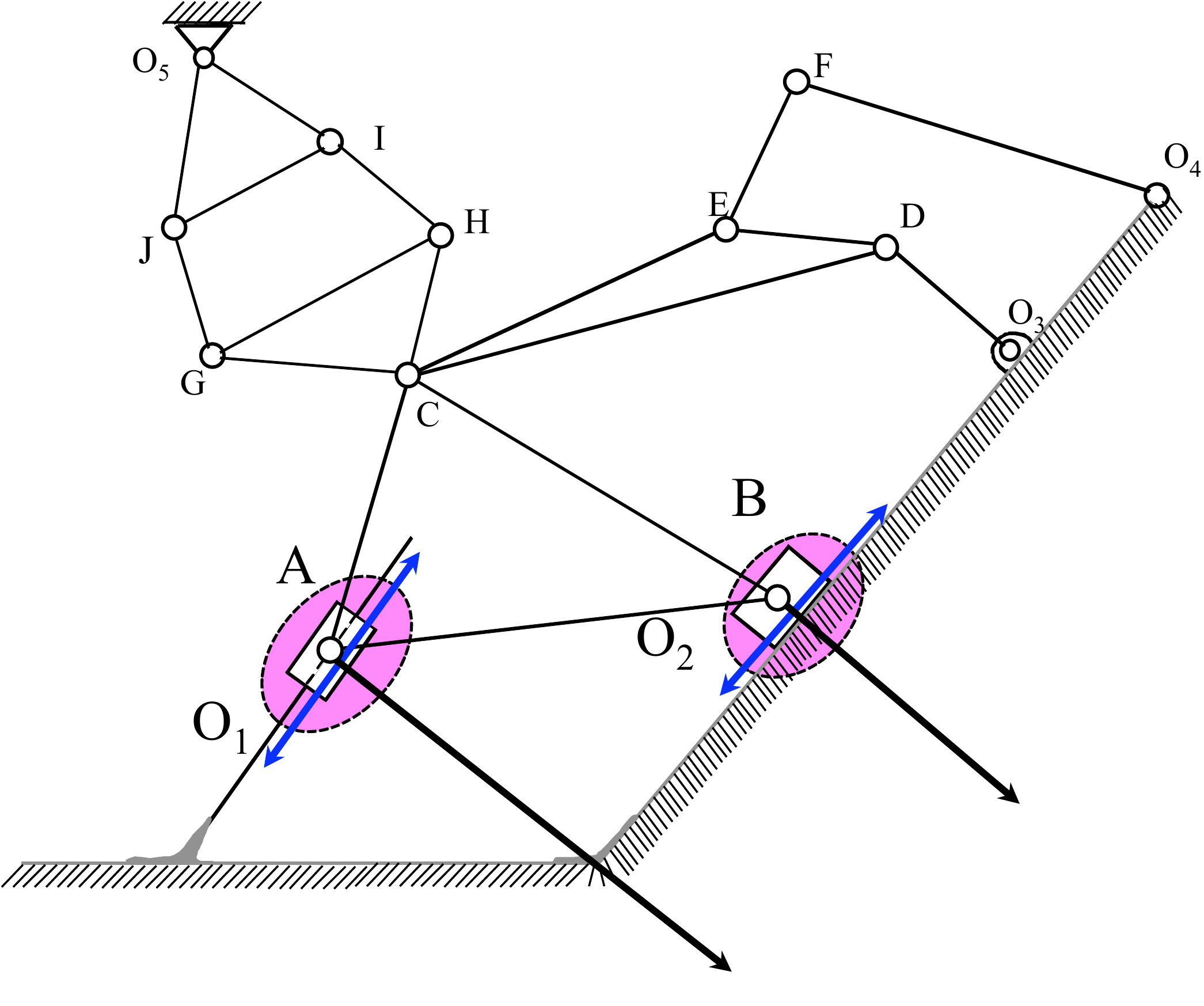}} \
  \subfigure[] {\includegraphics [width=.295 \textwidth]{StructuralScheme.pdf}}
 \caption{Slider joints in a linkage (a), function like joints to  infinite bars perpendicular to the motion (b) and are then translated to schematic bars and joints in the structural scheme (c).\label{fig:slide}}
\end{figure}

It is appropriate in an applied mathematics to say a bit
more about how some unusual features in the linkage were
translated in the graph and about framework constraints. Most links in
the plane linkage were represented as bars with pin joints at
their ends. The translation for these is clear:  links go to fixed
length bars and pin-joints go to vertices.  There were two more
exceptional cases that we need to address - highlighted in
Figure~\ref{fig:slide}.

In this transfer, the slider joints are first seen
as equivalent to a link to a pinned vertex very far away (at
`infinity'), which leaves only a translation (perpendicular to the
bar) (b) and then schematically represented as a pin at a finite
location (c).  Geometrically, the inclusion of slide constraints
as `projective points at infinity' is part of the mechanical
folklore and the corresponding projective theory of rigidity
\cite{CrapoWhiteley}.   In this literature, it is recognized that
the infinitesimal rigidity of a framework,  the rank of the
rigidity matrices, and the dimension of the space of non-trivial
velocities is projectively invariant
\cite{CrapoWhiteley,WChapter}.

It is also well known in the engineering literature
(\cite{Norton}) that both sliders and pinned joints are of the
same type - lower kinematic pairs. Therefore, from the view point
of rigidity (counts of the degrees of freedom) the sliders can be
treated as pinned joints as is done in the structural scheme.

With a finite framework, it is possible to use such a projective
transformation to bring all joints to finite locations.  The rest
of the schematic, and the analysis in \S2,\S3 is combinatorial, not
specifically geometric.   In the pinned rigidity matrices above,
and our work below, we will continue to use the `simpler'
Euclidean representation with finite points.  However,  there is a
full matrix representation and associated analysis that explicitly
includes joints at infinity.

\subsection{Motions generated by removing an edge: \\ \dassur vs \sdassur graphs}
The following result generates a 1 DOF linkage by removing an
edge from an $d$-isostatic framework.  Recall that an edge is part of a unique \dassur
graph in the extended decomposition, in which each strongly
connected component is extended to include its outward directed
edges. Recall that a {\em \sdassur  graph} is a \dassur  graph
(minimal pinned $d$-isostatic graph) with the added property that
removal of any edge induces a motion of all the inner vertices. In the plane all
$2$-Assur graphs are strongly $2$-Assur graphs, as we will show.

\begin{figure}[htb]
\centering
 \subfigure[] {\includegraphics [width=.45\textwidth]{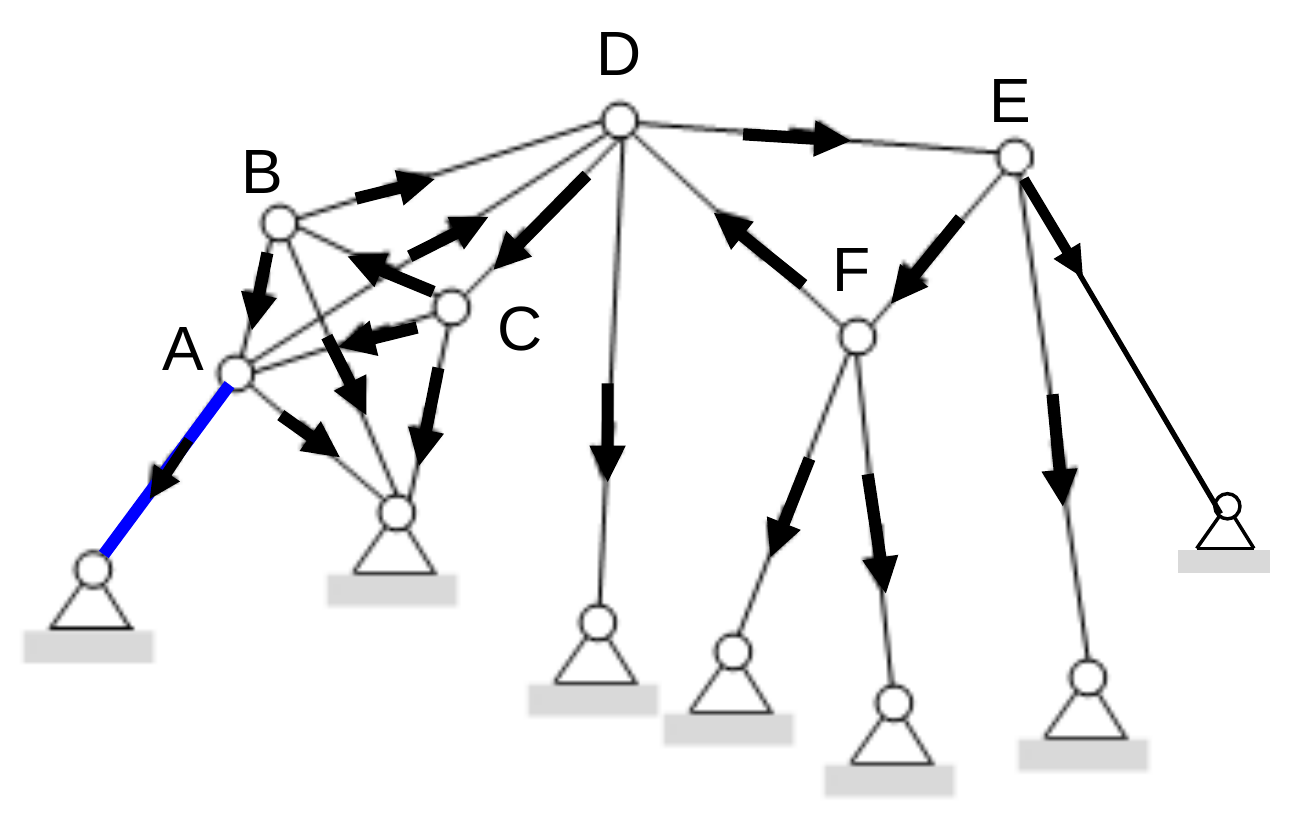}} \
  \subfigure[] {\includegraphics [width=.5
  \textwidth]{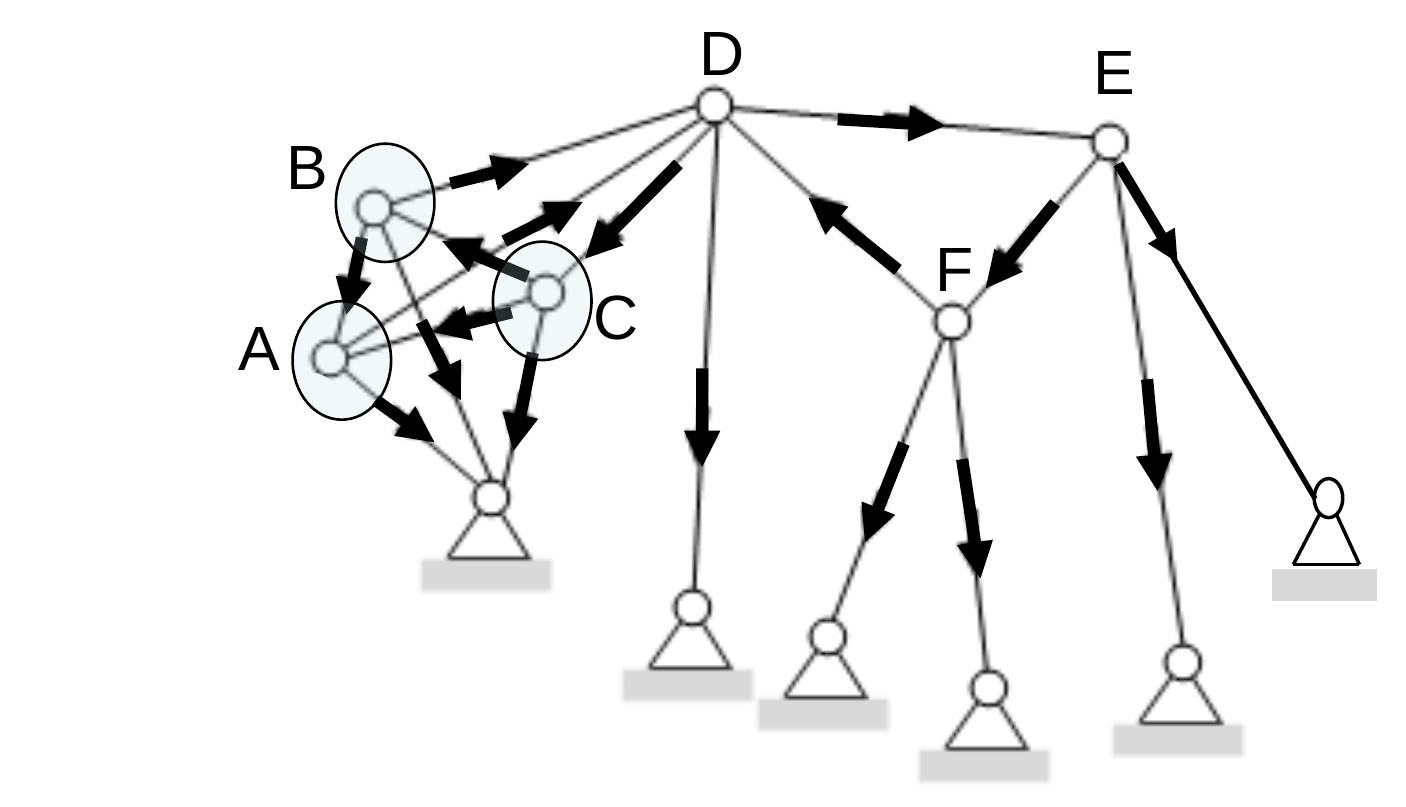}}\
  \subfigure[] {\includegraphics [width=.24
  \textwidth]{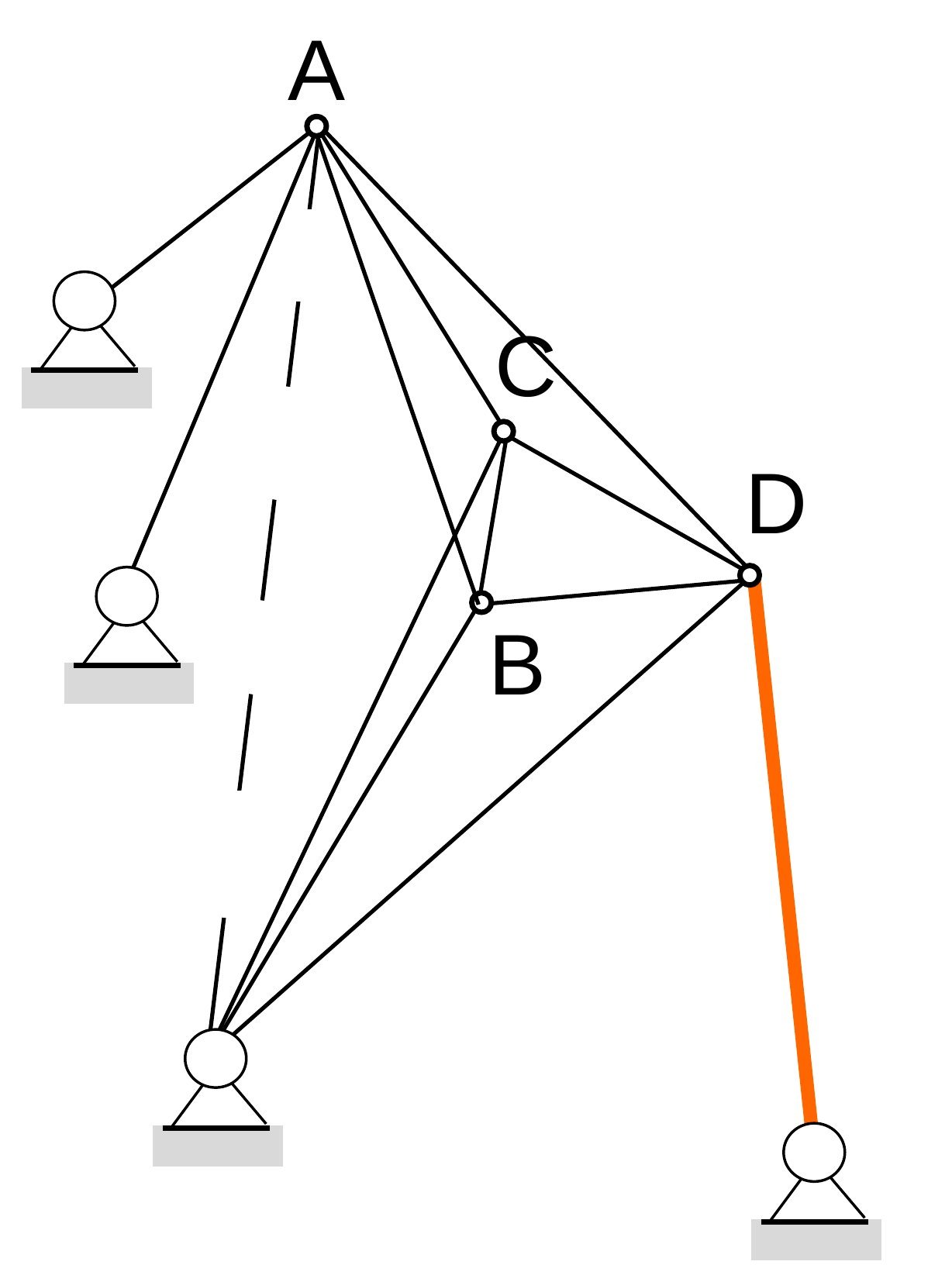}}\
    \subfigure[] {\includegraphics [width=.32
  \textwidth]{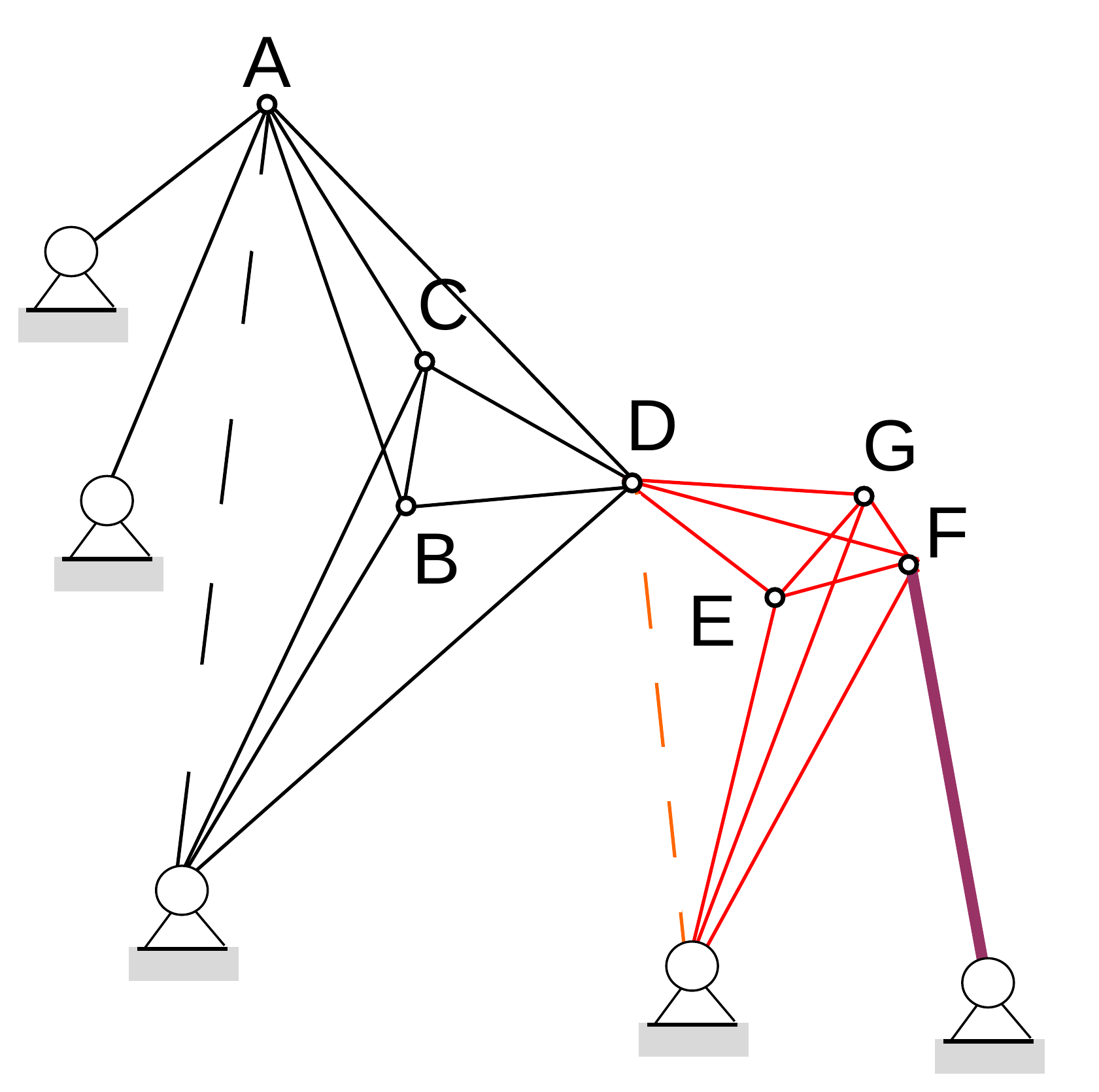}}\
    \subfigure[] {\includegraphics [width=.42
  \textwidth]{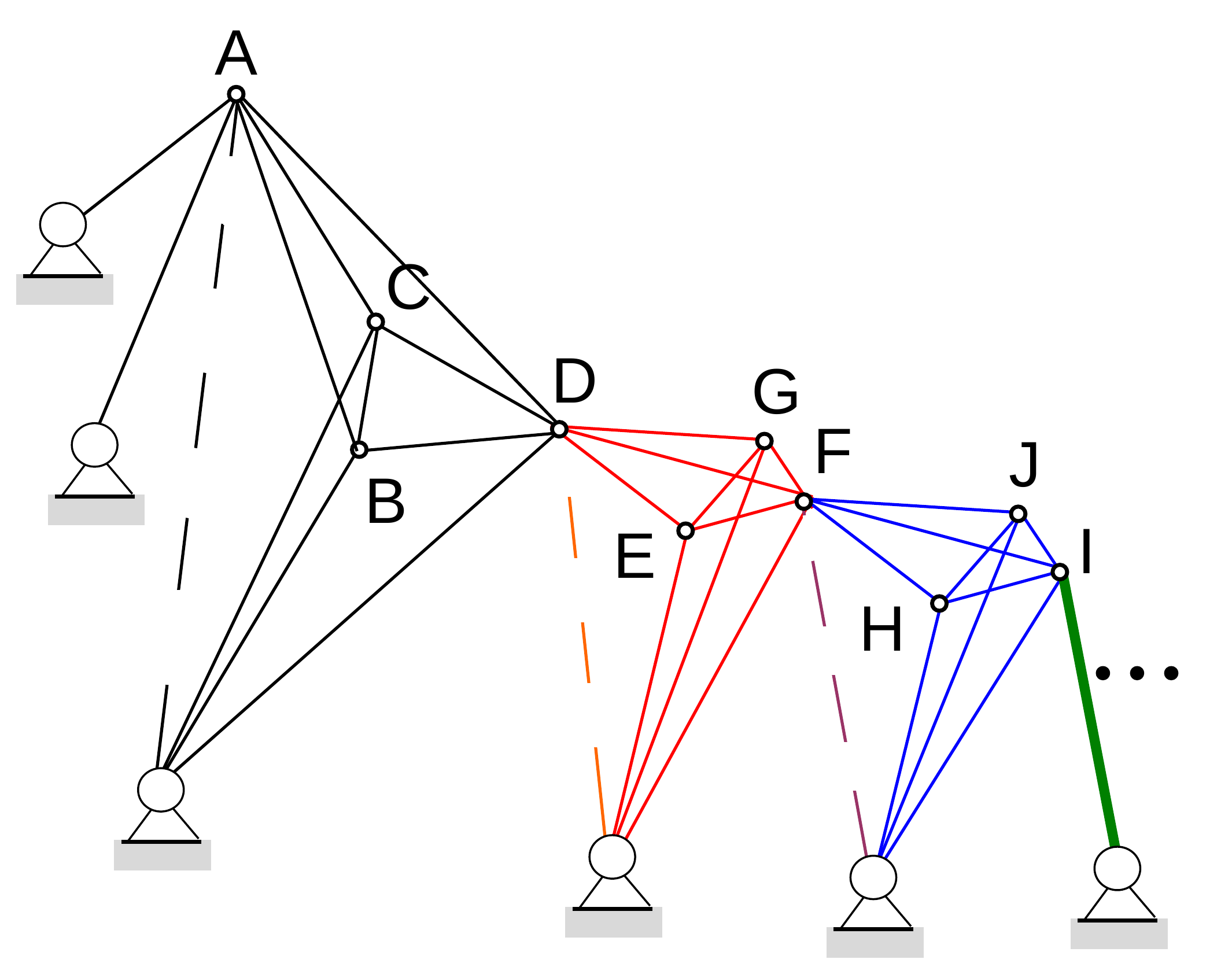}}
\caption{Examples of \thassur graphs which are not \sthassur.
Removal of a blue edge (driver) (a) results in a motion of only
the circled vertices in (b). In the examples in (c, d, e),
introduce a new copy of a banana graph, illustrate how removal of
certain edges induces smaller sets of vertices to move. Removal of
an orange edge in (c) causes all vertices to move except A, which
is held rigidity to the ground by an implicit (imaginary) edge
indicated with a dashed line. This process is continued through
(d,e). \label{fig:weakassur}}
\end{figure}

\begin{prop} If $\widetilde{G}$ is a pinned $d$-isostatic graph, and at a generic configuration $p$, removal of any edge from
$\widetilde{G}$ causes an infinitesimal motion which is non-zero
at all its inner vertices, then $\widetilde{G}$ is \sdassurp.
\label{prop:remove}
\end{prop}

\proof We need to show that $\widetilde{G}$ is a `minimal' pinned
$d$-isostatic graph. By assumption $\widetilde{G}$ is a pinned
$d$-isostatic graph. Assume $\widetilde{G}$ is not minimal pinned
$d$-isostatic. Then by Theorem \ref{thm:2blockdecomp} there is a
\blocktriangular matrix decomposition with more than one block.
Removing an edge from any block which is lower right in the matrix will leave
the graph associated with the upper left block (equivalently at the bottom of
the directed graph decomposition) as pinned $d$-isostatic.  This
guarantees that the solution $R(\widetilde{G},p) \times U = 0$ is
zero on all vertices of this upper block (i.e. these vertices have
no motion), a contradiction. Therefore $\widetilde{G}$ must be a
minimal pinned $d$-isostatic graph. Since $\widetilde{G}$ is
minimal and by assumption removal of any edge from $\widetilde{G}$
causes all inner vertices to go into motion, $\widetilde{G}$ is
\sdassur. \eop

Note that in Proposition \ref{prop:remove}, we could have assumed
that $\widetilde{G}$ is pinned $d$-rigid, and the same conclusion
would still follow. The assumption of  independence was not necessary, as
removal of an edge which causes a motion in the
graph, indicates that the edge is independent.

\begin{figure}[htb]
\centering {\includegraphics
[width=.60\textwidth]{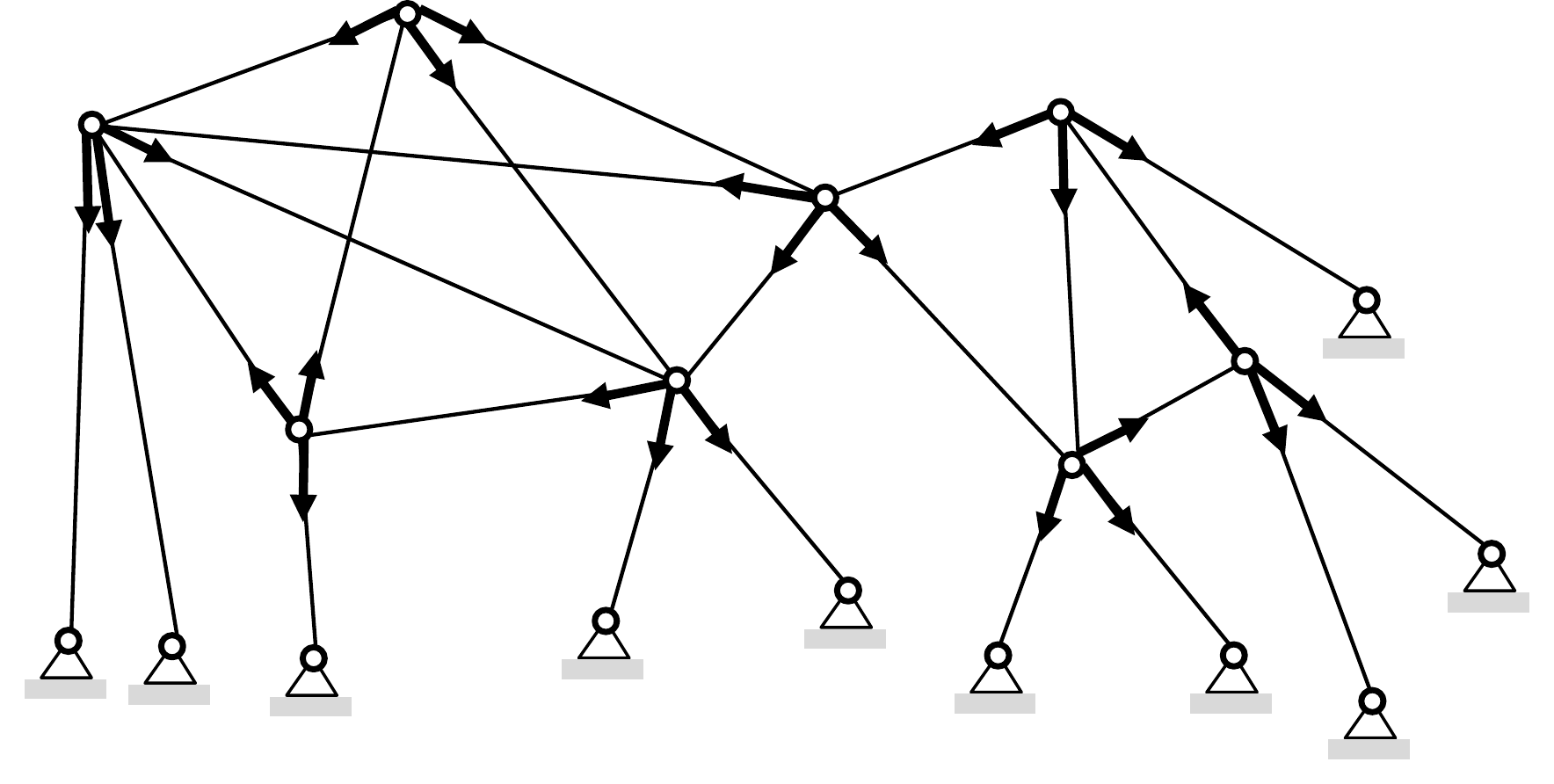}}
\caption{An example of a \sthassur graph. Removal of any ground edge
or inner edge 
induces a motion of all inner vertices.\label{fig:3dAssurgraphb}}
\end{figure}

We have the following stronger statement about  $2$-Assur
graphs (Figure~\ref{fig:assurweaklyassur}(b):

\begin{prop} If we remove any edge from a  $2$-Assur graph $\widetilde{G}$ then this leaves a graph which, at a generic configuration $p$ has an infinitesimal motion which is non-zero at all inner
vertices.  Therefore $\widetilde{G}$ is strongly $2$-Assur.  \label{prop:remove2D}
\end{prop}

\proof  Assume the graph $\widetilde{G}$ is  $2$-Assur. This
means that $|E|=2|I|$.  Removing any one edge will leave the
pinned rigidity matrix with $|E|=2|I| -1$, so there must be a
nontrivial solution $U$ to the matrix equation $R(\widetilde{G},p)
\times U = 0$. If $U_{i} = 0 $ on some inner vertex $i$, then that
vertex is still rigidly connected to the ground, and therefore
must be in a pinned subgraph $\widetilde{G}'$ with $|E'|=2|I'|$.  This
subgraph $ \widetilde{G}'$ would itself be a pinned $2$-isostatic
graph which could not include at least one of the vertices of the
removed edge.  Such a subgraph contradicts the minimality of the
original  $2$-Assur graph $\widetilde{G}$. \eop

The analog of  Proposition~\ref{prop:remove2D} in $d$ space ($d$
$>$ $2$) fails, and there are explicit counter-examples. The
example in Figure~\ref{fig:3dAssurgraphb} is a \sthassur graph
since removal of any edge causes all inner vertices to be mobile,
while the examples in Figure \ref{fig:weakassur} are only
\thassurp. We should point out that all the examples that were
presented in \S3 were of strongly 2 or 3-Assur graphs.

The examples in Figure \ref{fig:weakassur} (c, d, e) are
particularly interesting. Starting with a  3-Assur graph in
(c) we keep constructing new  3-Assur graphs as shown in
(d) and (e). In these set of examples, removal of some edges
(drivers) induces a full motion of inner vertices while removal of
other edges induces motions in a smaller collection of inner
vertices. One could create a further partial order of edges (drivers) in the
\thassur graph - those edges whose removal puts all the
inner vertices in motion which we call \emph{regular} drivers, and
further classification of edges called \emph{weak} drivers. The
partial order among the weak drivers would depend on the partial
order of subsets of inner vertices which are sent into motion by a
removal of this weak driver.

In 2D this
distinction among the drivers would only occur in the graphs that
are not 2-Assur, where removal of any edge in one component will
only cause the motion in that 2-Assur component and all the other
2-Assur components below it in the acyclic 2-Assur decomposition.

The fact that not all \thassur graphs are \sthassur is
connected to combinatorial obstacles to a good combinatorial
(counting) characterization for $3$-dimensional bar and joint
frameworks, also noted in the lack of  necessary and sufficient
counting conditions in higher dimensional frameworks. In this
context, detecting the difference between \dassur graphs
and \sdassur graphs can be equally challenging. We pose this
difficult open problem in combinatorially rigidity: {\em Given a
\dassur graph ($d$ $>$ 2) $\widetilde{G}$, find a combinatorial
method (using only the information from the graph) that determines
if $\widetilde{G}$ is \sdassur}.

\medskip
\noindent {\bf Remark.} Despite these difficulties in higher
dimensional space, having a comprehensive understanding of
2-dimensional Assur graphs can be very useful to the analysis of
3-dimensional linkages. A frequent practice by mechanical
engineers and in robotics community is to decompose the
3-dimensional structure to several 2-dimensional components, which
significantly simplifies the analysis.

In many cases,
3-dimensional linkages are built up by carefully connected copies of planar
structures. Often identical 2-dimensional structures are reused
forming the larger 3-dimensional structure, where the motions of
individual 2-dimensional structures are restricted to the plane by the geometry of the other constraints.
\eop

By general results from algebraic geometry  at a generic configuration $p$ for the
vertices (the length constraints
are algebraic conditions),  the infinitesimal motion $U$ extends to a non-trivial
finite path $p(t)$ within the configuration space $\R^{d|I|}$
which preserves the constrained lengths \cite{asiroth}.  In fact,
the mechanism is typically designed to create a specific finite
path in the configuration space, or perhaps for a given inner
vertex to trace a specific path in $\R^{d}$.
Thus, we can summarize the two propositions from this section in
terms of finite motions in the following corollaries:

\begin{cor} If $\widetilde{G}$ is a pinned $d$-isostatic graph, and at a generic configuration $p$, then removal of any edge from
$\widetilde{G}$ causes a finite motion which is non-zero at all
its inner vertices if any only if $\widetilde{G}$ is \sdassurp.
\label{thm:remove}
\end{cor}

One of the central motivations for decomposing a linkage is to
break down the analysis of the paths being generated from one set
of large polynomial constraints for the entire linkage into the
analysis of smaller polynomials for each of the pieces, plus a set of linking equations
for composing the results for the components into a single larger
analysis.  This is precisely what the \dassur
decomposition lets us do.  Here, we present a linearized version of this
process, for the decomposition into \dassur graphs.

Assume we have a $d$-isostatic linkage (pinned framework)
$(\widetilde{G_{k}},p)$, and we assign a vector $r$ of {\em drive
velocities} to the pinned vertices. We need the $e\times 1$ {\it
drive matrix}:
\begin{displaymath}
\mathbf{D}(\widetilde{G_{k}},p,r) = \bordermatrix{ &  \cr
& \vdots \cr
\{i,j\}\ |i,j\in I & 0  \cr
\{i,k\} \ |k\in P &  (p_{i}-p_{k})\cdot r_{k}  \cr
 &\vdots
 }
\end{displaymath}

\begin{prop} For an assignment of drive velocities $r$ to the pinned
vertices of a \dassur graph $\widetilde{G_{k}}$, and a
generic configuration $p$ for the vertices,  the drive equation
$\mathbf{R} (\widetilde{G_{k}},p)\times U =
\mathbf{D}(\widetilde{G_{k}},p,r)$ has a unique solution.
\end{prop}

\proof  The essential property is that  $\mathbf{R}
(\widetilde{G_{k}},p)$ is invertible, so the system of
non-homogeneous equations has the unique solution:
$$\quad\quad\quad\quad\quad\quad\hfill U=[ \mathbf{R} (\widetilde{G_{k}},p)]^{-1}     \mathbf{D}(\widetilde{G_{k}},p,r). \hfil \quad\quad\quad \quad \quad\quad\quad\quad\quad\quad\square$$

With this observation in hand, we see that if we decompose the
original $3$-isostatic linkage, and replace one edge of a bottom
\thassur linkage $(A_{1},p)$ with a driver, then we can
compute the velocities of all the inner vertices of $(A_{1},p)$.
These velocities can then be used as drive velocities for other
components in the linkage.  Iteration up the decomposition gives
solutions for all inner vertices of the entire linkage.  All
computations are reduced to computations within Assur components.
This is one primary value for the Assur
decompositions of linkages \cite{Assur,Mitsi,Norton}.

Notice that this capacity to propagate drive velocities through the decomposition does not require that all the components are \sdassurp.
However, it is possible that if some are not, the induced velocity will be zero at a number of inner vertices of components above the driver.

\medskip

\noindent {\bf Remark.}  If we go to a non-generic configuration
$p$ for a \dassur graph $A$, then the rank of the pinned
rigidity matrix can drop, creating a {\it singular configuration}.  In
these singular configurations, two things happen:   (i) there are drive
velocities for which there are no solutions; and (ii) for some
drive velocities (including the $0$ drive velocity), there are
multiple solutions.  Both of these events are a serious problem
for a linkage in mechanical engineering.

The configurations $p$ which make this happen can be {\it dead
end positions}, a geometric subclass of the singular positions - depending on the driver, or drive velocities.
There is an extensive geometric literature for dead end positions
\cite{Penne,SSWII}.  Again, this geometric analysis is simplified
by working with the geometry of the components of the Assur
decomposition.   We note that the singular positions of an Assur
component depend on the geometry of the configuration $p$ of both
the inner and the pinned vertices.   This geometry is the subject
for ongoing investigations, often for specific linkages or classes
of linkages.          \eop

\subsection{Re-pinning of inner vertices and release of pinned vertices}
\begin{figure}[htb]
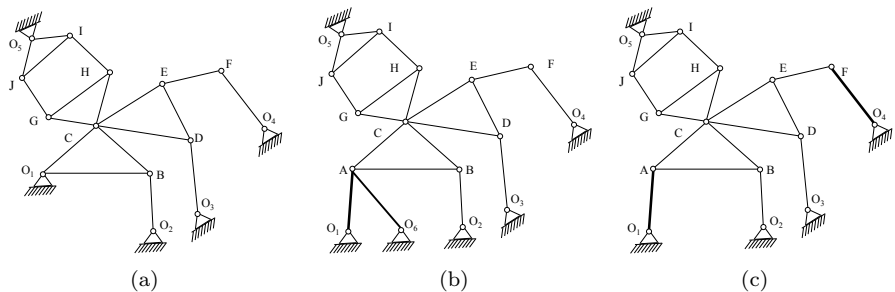

\centering
  \subfigure[] {\includegraphics [width=.32\textwidth]{IsostaticExample1.pdf}} \
  \subfigure[] {\includegraphics [width=.32 \textwidth]{AdaptedScheme1.pdf}}\
  \subfigure[] {\includegraphics [width=.32 \textwidth]{StructuralScheme.pdf}}
 \caption{Shifting from a pinned vertex $O_{3}$ (a), through a dyad at A (b), to a flexible structural scheme by removing an edge (c). \label{fig:liftpin}}
\end{figure}

An alternative operation which  mechanical engineers use in 2D
linkages is to `replace a driver' in a structural scheme and
create a 2-pinned isostatic graph is to shift an inner vertex into
a pinned vertex (Figure~\ref{fig:liftpin}(c) to  (a)). This can be
represented by a two step process passing through
Figure~\ref{fig:liftpin}(b).  We add a bar to freeze out the
motion from the inner end of the driver.  Then we convert this
small dyad (which includes the old driver) into a new pinned
vertex Figure~\ref{fig:liftpin}(a).

The reverse operation, which also applies in higher dimensional space,
is to `release' a pinned vertex into a new inner vertex, moving
from Figure~\ref{fig:liftpin}(a) to (b), using $d$ edges in
$d$-space ($d=2,3$) and creating a further \sdassur graph
(dyad in the plane, triplet in $3$-space). If we started with a
pinned $d$-isostatic framework, this released framework is also a pinned
$d$-isostatic framework, with one added \sdassur
component (b).   Then we remove
an edge from this minimal \dassur graph (dyad) (c),
generating a local motion at this new inner vertex. This motion
propagates on into rest of the inner vertices, using
the matrix as above.

More generally, we can take a $d$-isostatic graph, and add a
single velocity to any one pin (or set of pins).  This then
proceeds as above with the drive matrix connecting these pins to
the other parts of the framework, and solving for the velocities
at all the inner vertices.  In following this propagation, it is still
valuable to be able to do the work on one \dassur
component at a time, which was the essential feature of the
previous subsection.

\subsection{Vertex removal}
An additional  engineering technique in testing for decompositions
is to remove an inner vertex. Removing any pinned vertex will
generate a finite space of non-trivial motions - the dimension of
the space being the valence of the pinned vertex. We now focus on removing an inner vertex. Any edge will always
contain at least one inner vertex, leading to a vertex analog of
Theorem~\ref{thm:remove}.

\begin{theorem} If $\widetilde{G}$ is a pinned $d$-isostatic graph, and at a generic configuration $p$, removal of any inner
vertex from $\widetilde{G}$ causes a finite motion which is non
zero at all inner vertices, then $\widetilde{G}$ is
\dassurp.\label{thm:removevertex}
\end{theorem}

\proof Case 1: $\widetilde{G}$ has one inner vertex.
$\widetilde{G}$ is trivially   \sdassurp.

Case 2: There is more than one inner vertex. We want to show that
$\widetilde{G}$ is minimal pinned $d$-isostatic. If some vertex has valence $d$, removing
it will not cause the other vertices to move, which contradicts
the assumption.  Therefore the vertex we remove from $\widetilde{G}$ has valence at
least $d + 1$.

Since $\widetilde{G}$ is a pinned $d$-isostatic graph, there is a
\blocktriangular matrix decomposition. Assume $\widetilde{G}$ is
not minimal pinned $d$-isostatic. Removing any vertex of degree
$\geq d+1$ from any block which is lower right, will leave the upper
left block (near the bottom of the directed graph decomposition)
as pinned $d$-isostatic. This guarantees that the solution
$R(\widetilde{G},p) \times U = 0$ is zero on all vertices of this
upper block (i.e. these vertices have no motion), a contradiction.
Therefore $\widetilde{G}$ is minimal and \dassurp.\eop.

In Theorem \ref{thm:removevertex} it may be surprising that we cannot
conclude that the graph will be \sdassur ($d$ $>$ 2). Consider
the example in Figure \ref{fig:weakassur} (a), removing any inner
vertex from this graph will cause all other inner vertices to be
mobile, yet this is a \thassur graph, not \sthassurp.

\begin{theorem}
If we remove any inner vertex from a \sdassur graph
$\widetilde{G}$, then at a generic configuration $p$,
$\widetilde{G}$ has a finite motion which is non zero at all inner
vertices. \label{thm:removevertexdassur}
\end{theorem}

\proof Case 1: $\widetilde{G}$ has one inner vertex, removing it
leaves no inner vertices.

Case 2: There is more than one inner vertex.  Since $\widetilde{G}$ is
\sdassur, every inner vertex has to be at least valence $d+1$,
as any inner vertex of degree $d$ (and its outgoing edges) in a
pinned $d$-isostatic graph is a \sdassur component (a $d$-dyad)).
Choose any inner vertex $v$ and remove its edges. As
$\widetilde{G}$ is \sdassurp, removal of these edges, in fact any
single edge, will cause a motion at all inner vertices. Now remove
the vertex $v$, all other inner vertices are still in motion. \eop

This result is expected, as removal of an inner vertex
should cause at least as much motion as a removal of a single edge
incident to that vertex. In particular, if we remove a vertex of
degree higher than $d+1$ there will be more flexibility (more DOF)
caused then removal of any edge at that vertex.

\begin{figure}[htb]
\centering {\includegraphics
[width=.45\textwidth]{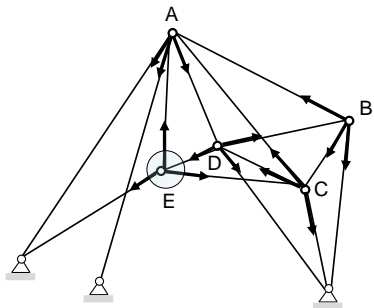}}
 \caption{A \thassur graph, removal of vertex E causes motion in vertices B, C and D, but not in A. \label{fig:removevertex}}
\end{figure}

\noindent {\bf Remark.} 
Most of these concepts and results in \S3,4 were illustrated using examples in 3-space.
All these observations are true in all dimensions. To give higher dimensional
examples, one can easily start with a 3-dimensional example and
construct similar examples in any higher dimension, using the
technique of coning \cite{WhiteleyConing}, which transfers the
rigidity, orientation (out-degree), counts, etc. from a framework
in dimension $d$ to dimension $d+1$.\eop

\section{Further Work}\label{sec:furtherwork}

The results in this paper are  contributions to the
the synthesis and analysis of linkages in $d$-space, and provide
new extensions of $2$-dimensional Assur graphs to higher
dimensional space, and new connections and better understanding of
the difficulties with higher dimensional rigidity. There are
several key further directions that are part of the ongoing
research, some of which are addressed in papers such as
\cite{mechanicalpebblegame} and others which hope to follow up in future
papers. We highlight some of current research and envision other
developments.

\subsection{Applications to more general kinematic structures}\label{subsec:body}
Section 3 developed decompositions for minimal pinned
isostatic frameworks (\dassurp), which traced back to the
strongly connected decompositions of Section 2.  Looking back, the
techniques as developed only depended on two properties:
\begin{enumerate}
\item we had an underlying constraint multi-graph; \item we have a
square constraint matrix for the multi-graph with rows indexed by
the edges and columns indexed (in groups) by the vertices.
\end{enumerate}
With these two properties, we can use the constraint matrix to
generate directions on the multi-graphs, with the out-degree
corresponding to the number of columns for a given vertex (as in
\S3).   With this in place, the Laplace block decomposition of the
determinant of the constraint matrix facilitates the certificate
that we can generate an orientation of the constraint multi-graph.
This orientation will yield a strongly connected decomposition of
the original multi-graph.

The key results in \S2 and \S3 also did not depend on each type of
vertex having the same degree in the constraint graph
(equivalently, having the same number of columns in the pinned
constraint matrix).  They only depended on the out-degree being
constant when we switched to another orientation of the graph or
equivalently on there being an associated number of columns or
variables for each type of vertex in the graph.  The results also
extend to multi-graphs as the constraint input, as long as there
is a row in the constraint matrix for each edge of the
multi-graph.   Everything generalizes directly to these broader
settings.

We can apply the entire suite of decomposition techniques from
general pinned bar and joint frameworks to analyze isostatic
body-bar frameworks in $3$-space\cite{WhiteWhiteleyII}. For these special frameworks, we
not only have the required constraint matrix on the multi-graph,
but we have complete necessary and sufficient counting properties
to fully test the components of the decomposition for whether they
are generically $3$-isostatic \cite{WhiteWhiteleyII}.  In fact,
there is a complete theory, with good algorithms, and these are
associated with some types of linkages, such as the well-known
Stewart platform.  In addition, for $d$-body-Assur graphs, all minimal components are strongly $d$-Assur, in the sense that removal of any one edge causes all inner bodies to go into relative motion.

\begin{figure}[htb]
\centering
   \subfigure[] {\includegraphics [width=.45\textwidth]{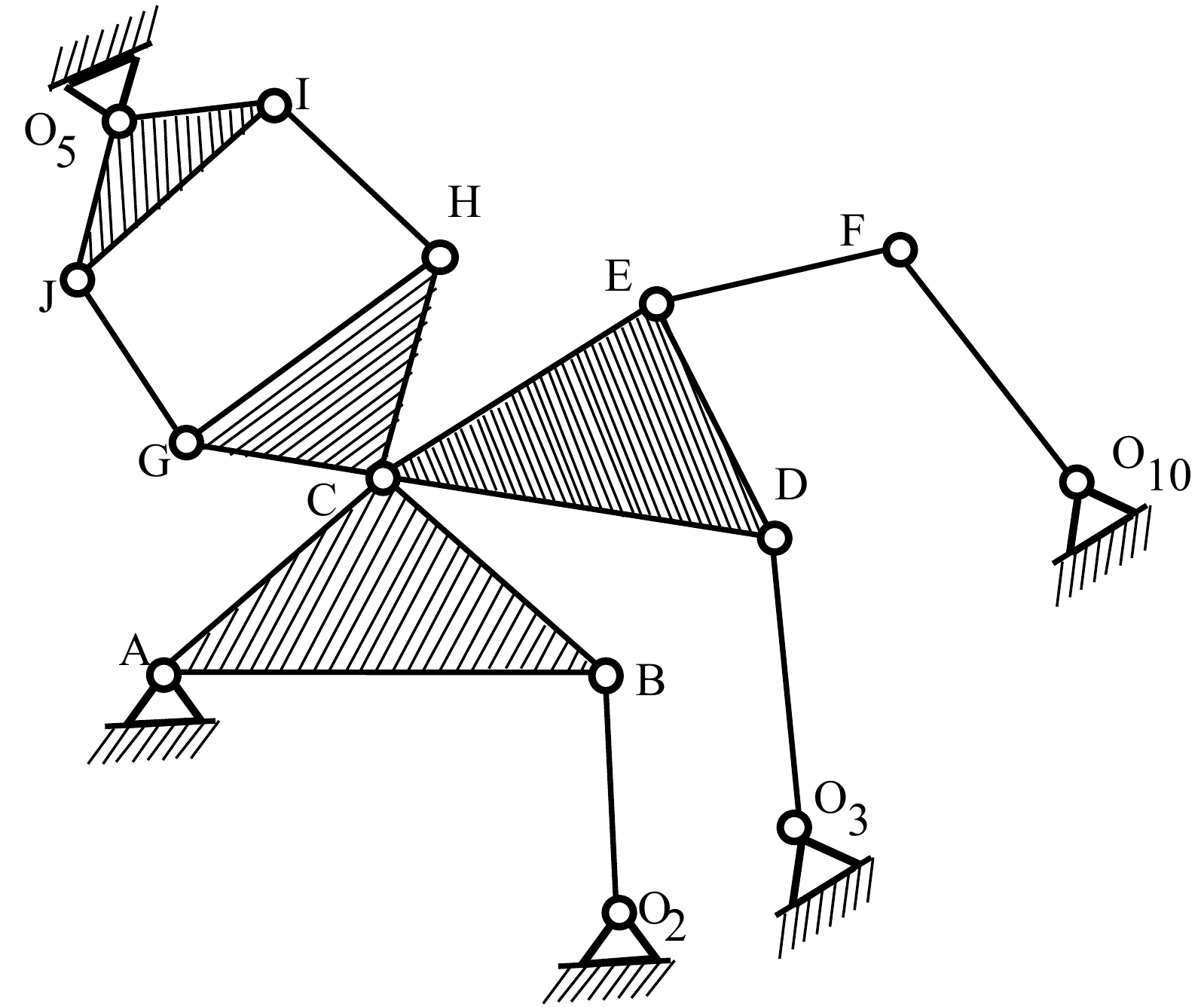}}\quad
 \subfigure[] {\includegraphics [width=.45\textwidth]{IsostaticExample1.pdf}}
\caption{The unified structural scheme (b) of the linkage with
bodies (a).\label{fig:bodylinkage}}
\end{figure}

We could also adapt these decomposition techniques to mixed
frameworks in the plane, in which we have bodies, bars, joints,
pins between bodies, prismatic or pistons drivers
(Figure~\ref{fig:bodylinkage}).   For these more general kinematic
chains, the classical counting rules are associated with the name
of Gr\"ubler \cite{Grubler,SSWI}.  These counting rules provide
necessary conditions on a linkage to be generically isostatic.
While they are not usually elaborated into necessary and
sufficient conditions in the engineering literature, they can be
reworked using the generic theory of plane rigidity to give such
necessary and sufficient conditions. All of the results for the
plane decomposition extend to mixed frameworks, providing a complete theory of
$2$-Assur decompositions, and the associated \blocktriangular
decomposition of the pinned rigidity matrix for the constraint
multi-graph.

We can  further  adapt these decomposition techniques to mixed
frameworks in $3$-space, with bodies, bars, joints, revolute
joints between bodies (hinges), and prismatic joints
\cite{Grubler,Norton}.  For these general kinematic multi-loop
mechanisms, the more general counting rules are associated with
the names of Chebychev, Gr\"ubler and Kutzbach \cite{Grubler}.  In
such a broad extension, the prismatic joints again become
`revolute joints at infinity' in the general projective theory of
$3$-space rigidity \cite{CrapoWhiteley}.

\subsection{Algorithms for generating directions and Assur decompositions}
The algorithms for assigning directions to the graph presented in this paper are, in principle,
exponential: searching over terms in the Laplace expansion of the
determinant of the pinned rigidity matrix. However, fast pebble
game algorithm can be used to generate the $2$-directed
orientations, as well as for the necessary and sufficient counting
conditions and to decompose the graph into 2-Assur components.
Some of this work using the pebble game algorithm and techniques
from rigidity theory was for the first time presented to a
mechanical engineering community \cite{mechanicalpebblegame},
which provides fast analysis and synthesis of linkages.

As we have seen, for $d \geq 3$-space, there is no known necessary
and sufficient counting (polynomial time) algorithm for infinitesimal rigidity or independence. The examples
in \S\ref{subsec:counts} illustrated some of the problems.
However, if we start with a pinned $d$-isostatic graph (or at
least without redundance - no dependant rows in the pinned
rigidity matrix), we have shown elsewhere how to adapt the pebble game to
give the $d$-directions and \dassur decomposition
\cite{mechanicalpebblegame}. We are currently investigating more
on these pebble games techniques.

For general body-bar frameworks in $d$-space, as well
as their specialization as body-hinge and molecular frameworks,
mentioned in \S\ref{subsec:body}, there are full,  efficient
pebble games to generate the directed graphs and to complete the
decomposition.
Body and bar frameworks are routinely used in study of general
linkages and in robotics. Given this motivation of important
applications to mechanical engineering community, we are currently
adapting some of the techniques and algorithms to body bar Assur
graphs.

\subsection{Concluding remarks}

As this paper confirms, there are a wide variety of important, and
mathematically interesting unsolved questions for investigation.
In the last 30 years, there has been a broad development in the
general theory of rigid and flexible structures. In particular,
during the last decade, there has been a renewed interaction of
these developments in the mathematical theory of rigidity, and the
results on rigidity of structures, with the parallel analysis of
linkages, which has its own rich history.  Our work on Assur graphs,
offers some additional tools for decomposition of
pinned bar and joint structures,  probing into the difficult and not well
understood rigidity of bar and joint structures in higher
dimensional space.

There is a wide field of fruitful directions for further
investigation.  We invite the reader to join in these
investigations.


\end{document}